\documentclass[leqno,12pt]{article}  
\usepackage{amsmath}
\usepackage{amssymb}

\usepackage[german,english]{babel}
\usepackage{amsthm}

\title{The \v{C}ech Filtration and Monodromy in Log Crystalline Cohomology}

\author{\textsc{Elmar Grosse-Kl\"onne}}
\date{}

\theoremstyle{plain} 
\newtheorem{satz}{Theorem}[section]  
\newtheorem{kor}[satz]{Corollary}  
\newtheorem{lem}[satz]{Lemma}  
\newtheorem{pro}[satz]{Proposition}  
\newcommand{\ho}{\mbox{\rm Hom}}  
 
\newcommand{\spec}{\mbox{\rm Spec}}  

\newcommand{\quot}{\mbox{\rm Quot}}  
\newcommand{\spf}{\mbox{\rm Spf}}  

\newcommand{\ord}{\mbox{\rm ord}}  
\newcommand{\bi}{\mbox{\rm Im}}  
\newcommand{\ke}{\mbox{\rm Ker}}  
\newcommand{\koke}{\mbox{\rm Coker}}  

\newcommand{\dlog}{\mbox{\rm dlog}}
\newcommand{\Gr}{\mbox{\rm Gr}}

\theoremstyle{remark}

\theoremstyle{definition}

\begin{document}
\maketitle
\footnote[0]{Most of this work was done during my visit at the University of California, Berkeley. I wish to thank Robert Coleman (and Bishop) for welcoming me there so warmly. Thanks are also due to Ehud de Shalit, Yukiyoshi Nakkajima and Arthur Ogus for useful related discussions. I thank the referee for his careful reading of the manuscript and his suggestions for improving the exposition. I am grateful to the Deutsche Forschungsgemeinschaft for supporting my stay at Berkeley.}

\begin{abstract}
For a strictly semistable log scheme $Y$ over a perfect field $k$ of characteristic $p$ we investigate the canonical \v{C}ech spectral sequence $(C)_T$ abutting to the Hyodo-Kato (log crystalline) cohomology $H_{crys}^*(Y/T)_{\mathbb{Q}}$ of $Y$ and beginning with the log convergent cohomology of its various component intersections $Y^i$. We compare the filtration on $H_{crys}^*(Y/T)_{\mathbb{Q}}$ arising from $(C)_T$ with the monodromy operator $N$ on $H_{crys}^*(Y/T)_{\mathbb{Q}}$. We also express $N$ through residue maps and study relations with singular cohomology. If $Y$ lifts to a proper strictly semistable (formal) scheme $X$ over a finite totally ramified extension of $W(k)$, with generic fibre $X_K$, we obtain results on how the simplicial structure of $X_K^{an}$ (as analytic space) is reflected in $H_{dR}^*(X_K)=H_{dR}^*(X_K^{an})$.
\end{abstract}

\begin{center} {\bf Introduction}\end{center}
Let $A$ be a complete discrete valuation ring of mixed characteristic $(0,p)$, with perfect residue field $k$ and fraction field $K$, and let $K_0$ be the fraction field of the ring of Witt vectors $W(k)$ of $k$. Let $X$ be a proper strictly semistable $A$-scheme. Besides its Hodge filtration the de Rham cohomology $H_{dR}^*(X_K)$ of the generic fibre $X_K$ of $X$ comes with a $K_0$-lattice with a Frobenius operator $F$ and a nilpotent operator $N$: these are obtained, via the Hyodo-Kato isomorphism (which depends on the choice of a uniformizer $\pi$ in $A$)\begin{gather}H_{dR}^*(X_K)\cong H^*_{crys}(Y/T)_{\mathbb{Q}}\otimes_{K_0}K\tag{$*$}\end{gather}from the Frobenius operator $F$ and the nilpotent operator $N$ on the Hyodo-Kato (log crystalline) cohomology $H_{crys}^*(Y/T)_{\mathbb{Q}}$ of the special fibre $Y$ of $X$. It follows from the theorem of Tsuji and (independently) Faltings that $H_{dR}^*(X_K)$, together with its Hodge filtration and the operators $F$ and $N$ on its $K_0$-lattice defined by $(*)$, allows the reconstruction of the $p$-adic \'{e}tale cohomology group $H_{et}^*(X_{\bar{K}},\mathbb{Q}_p)$ together with its Gal$(\bar{K}/K)$-action. We have $N=0$ if $X$ has good reduction. It is a general and important problem to reconstruct as much as possible of the $N$-structure on $H_{dR}^*(X_K)=H_{dR}^*(X_K^{an})$ (where $X_K^{an}$ is the rigid analytic space associated with $X_K$) solely from $X_K$ or $X_K^{an}$. For this purpose the obvious idea is to look at the following spectral sequence $(C)_S$. Let $\{Y_j\}_{j\in R}$ be the set of irreducible components of $Y$. Each $Y_j$ is (classically) smooth over $k$. For $i\ge 1$ let $Y^i=\coprod_{|I|=i}(\cap Y_{j})_{j\in I}$ where $I$ runs through the set of subsets of $R$ with precisely $i$ elements. We assume that all connected components of $Y$ are of the same dimension $d$. For direct sums $E=\coprod E_s$ of subschemes $E_s\subset Y$ let $]E[_X=\coprod_s]E_s[_X$ be the direct sum of the preimages of the $E_s$ under the specialization map $X_K^{an}\to Y$: these $]E_s[_X$ are admissible open subspaces of $X_K^{an}$. The admissible open covering $X^{an}_K=\cup_{j\in R}]Y_j[_X$ is a covering by {\it contractible} spaces in the sense of Berkovich. It gives rise to the spectral sequence\begin{gather}E_{1}^{pq}=H^q_{dR}(]Y^{p+1}[_{X})\Longrightarrow H^{p+q}_{dR}({X}_K).\tag*{$(C)_S$}\end{gather} In applications, for example when $X_K$ is a Shimura variety, $(C)_S$ often has arithmetical meaning. In the case $d=1$, Coleman and Iovita \cite{coliov} gave a description of $N$ on $H^{1}_{dR}({X}_K)$ in terms of $(C)_S$; namely, they proved that it is the composite$$H^{1}_{dR}({X}_K)\stackrel{Res}{\longrightarrow}H_{dR}^0(]Y^2[_X)\stackrel{\delta}{\longrightarrow}H^{1}_{dR}({X}_K)$$where $Res$ is an explicit residue map and $\delta$ is the connecting homomorphism in $(C)_S$. For any $d$, Alon and de Shalit \cite{alsh} gave a tentative definition of $N$ for varieties $X_K$ uniformized by Drinfel'd's symmetric spaces; their central concept of harmonic cochains on the Bruhat Tits building is intimately related to $(C)_S$. The observation from \cite{hkstrat} that for such varieties the filtration $F_C^{\bullet}$ on $H^{d}_{dR}({X}_K)$ defined by $(C)_S$ is the weight filtration for the Frobenius action plays a role in de Shalit's recent proof of the monodromy weight conjecture for such varieties \cite{desh} (another proof was given by Ito \cite{ito}). Besides these example we are not aware of other investigations of $(C)_S$. The general reference for monodromy operators arising in semistable families (considered there in the $\ell$-adic and in the complex analytic setting) is Illusie's article \cite{ilau}. Specifically, the Cech complex from \cite{ilau} 2.1.5, 3.2 does {\it not} correspond to $(C)_S$ but rather to the canonical Cech spectral sequence\begin{gather}E_1^{pq}=H_{rig}^q(Y^{p+1})\Longrightarrow H_{rig}^{p+q}(Y)\tag*{$(C)_{rig}$}\end{gather}with $H_{rig}^*$ denoting non logarithmic rigid cohomology. 
 
In this paper we begin the study of how $N$ interacts with $(C)_S$ for general $X$ by using log convergent cohomology. Since $Y$ is a normal crossings divisor on the regular scheme $X$ it gives rise to a natural log structure on $X$, and by pull back to natural log structures on subschemes of $X$. Thus, since $Y$ is a log scheme over the log point $T_1=(\spec(k),(\mathbb{N}\to k, 1\mapsto0))$ over $k$, all subschemes of $Y$ become $T_1$-log schemes. As such we ask for their log convergent cohomology $H_{conv}^*(./T)$ relative to $T=(\spf(W(k)),(\mathbb{N}\to W(k), 1\mapsto0))$ which takes values in $K_0$-vector spaces. For example, $H^*_{conv}(Y/T)=H^*(Y,{\bf C}\omega_Y^{\bullet})$ for a certain logarithmic de Rham complex ${\bf C}\omega_Y^{\bullet}$ on analytic tubular neighbourhoods of $Y$ in local $T$-log smooth liftings. Since $Y$ is log smooth over $T_1$ the Hyodo-Kato cohomology $H_{crys}^*(Y/T)_{\mathbb{Q}}$ is isomorphic to $H^*_{conv}(Y/T)$. There is a natural \v{C}ech spectral sequence\begin{gather}E_{1}^{pq}=H^q_{conv}(Y^{p+1}/T)\Longrightarrow H^{p+q}_{conv}(Y/T)=H_{crys}^{p+q}(Y/T)_{\mathbb{Q}}.\tag*{$(C)_T$}\end{gather}By \cite{hkstrat}, 3.8 and 3.9, the isomorphism $(*)$ extends to an isomorphism of spectral sequences $(C)_S\cong (C)_T\otimes_{K_0}K$. In particular, the descending filtration $(F_C^rH_{crys}^*(Y/T)_{\mathbb{Q}})_{r\ge 0}$ on $H_{crys}^*(Y/T)_{\mathbb{Q}}$ induced by $(C)_T$ redefines the descending filtration $(F_C^rH_{dR}^*(X_K))_{r\ge 0}$ on $H_{dR}^*(X_K)$ induced by $(C)_S$ (we call them the {\it canonical \v{C}ech filtrations}). As our main tool to analyse $(C)_T$ we introduce the {\it \v{C}ech-double complex} $B^{\bullet\bullet}$ on $Y$ (or rather on a simplicial scheme $U_{\bullet}$ associated with a suitable Zariski open covering of $Y$) whose total complex $B^{\bullet}$ computes $H_{conv}^*(Y/T)$ and whose $(k-1)$-st column $B^{\bullet k-1}$ computes $H^q_{conv}(Y^{k}/T)$. (Since $Y^k/T_1$ is not log smooth it is doubtful if its log {\it crystalline} cohomology (relative to $T$) is a useful object.) This means that $(C)_{T}$ is the spectral sequence associated with the stupid vertical filtration $F_C^{\bullet}$ of $B^{\bullet\bullet}$. On the other hand we consider a double complex $A^{\bullet\bullet}$ with associated total complex $A^{\bullet}$ which is the straighforward analog (in the log convergent setting, as opposed to the log crystalline setting) of the Hyodo-Steenbrink-complex $W_{\bullet}A^{\bullet}$ on $Y$ constructed by Mokrane \cite{mokr}. This $A^{\bullet}$ computes $H_{crys}^*(Y/T)_{\mathbb Q}=H_{conv}^*(Y/T)$, too. Moreover $A^{\bullet}$ is endowed with an endomorphism $\nu$ which induces $N$ on $H_{conv}^*(Y/T)_{\mathbb{Q}}$. The complexes $A^{\bullet}$ and $B^{\bullet}$ can be related in two ways: Firstly, there is a natural quasiisomorphism $\psi:A^{\bullet}\to B^{\bullet}.$  Secondly, there is a third complex $C^{\bullet}$ endowed with a filtration $F_C^{\bullet}$ and an endomorphism $\nu$ such that: there is a natural quasiisomorphism $A^{\bullet}\to C^{\bullet}$ respecting the respective endomorphisms $\nu$, and there is a natural quasiisomorphism $B^{\bullet}\to C^{\bullet}$, filtered with respect to the respective filtrations $F_C^{\bullet}$. In particular, $\nu$ on $C^{\bullet}$ induces an operator $N$ on $(C)_{T}$. We prove (cf. \ref{mokomptri}, \ref{inkl}, \ref{viawb}):

{\bf Theorem:} {\it $N=0$ on $H^*(Y,B^{\bullet k-1})$ and $N(F_C^{k-1}H^*_{crys}(Y/T)_{\mathbb Q})\subset F_C^{k}H^*_{crys}(Y/T)_{\mathbb Q}$ for any $k\ge0$. For $i\ge 0$ the $i$-fold iterated monodromy operator $N^i$ on $H^*_{crys}(Y/T)_{\mathbb Q}=H^*(Y,{\bf C}\omega_Y^{\bullet})=H^*(Y,B^{\bullet})$ is induced by a composite map of sheaf complexes$${\bf C}\omega^{\bullet}_{Y}\stackrel{\rho_i}{\longrightarrow} F_C^iB^{\bullet}\subset B^{\bullet}$$with an explicit residue map $\rho_i$.}

Retransposing to our lifted situation we derive (cf. \ref{kinkl}):\\

{\bf Theorem:} {\it $N$ on $H^{*}_{dR}({X}_K)$ naturally extends to an operator $N$ on $(C)_S$. However, $N=0$ on $H^*_{dR}(]Y^{k}[_{X})$, and $N(F_C^{k-1}H^{*}_{dR}({X}_K))\subset F_C^k H^{*}_{dR}({X}_K))$ for any $k\ge1$. For $i\ge 0$ the $i$-fold iterated monodromy operator $N^i$ on $H^{*}_{dR}({X}_K)$ has the form $$H^{*}_{dR}({X}_K)\stackrel{Res}{\longrightarrow}F_C^i H^{*}_{dR}({X}_K)\subset H^{*}_{dR}(X_K)$$with a residue map $Res$.}\\

We view this as the generalization, to any $d$, of the description of $N$ given by Iovita and Coleman in case $d=1$. In particular we get an upper bound for the vanishing order of $N$ described in terms of $X_K^{an}$. We do not know if in general the residue map $H^{*}_{dR}({X}_K)\stackrel{Res}{\longrightarrow}F_C^i H^{*}_{dR}({X}_K))$ can be made explicit {\it without} involving the log basis $T$. However, for $i=*=d$ there is a natural candidate, generalizing the residue map used by Iovita and Coleman in case $i=*=d=1$.

The inclusion $\bi(N^k)\subset F_C^kH^{*}_{dR}({X}_K)$ is an equality if for example the canonical map $H_{rig}^*(Y)\to H_{crys}^*(Y/T)_{\mathbb{Q}}$ is {\it strict} with respect to canonical \v{C}ech filtrations, which on $H_{rig}^*(Y)$ is defined through the above spectral sequence $(C)_{rig}$ (see Proposition \ref{dannja}). In general, however, it is {\it not} an equality, even if $Y$ is projective, $k$ is finite and the monodromy weight conjecture holds, see section \ref{gleich} for a counter example.

Yet there is more structure on $B^{\bullet}$: if $k$ is finite and the monodromy weight conjecture holds, the monodromy filtration on $H_{crys}^*(Y/T)_{\mathbb{Q}}$ is induced from a filtration on $B^{\bullet}$ (Theorem \ref{fiwbmo}). Moreover $B^{\bullet}$ has a product structure inducing Poincar\'{e}-duality (in contrast to $A^{\bullet}$, it seems), see section \ref{doubledef}.

We also ask how $N$ is related to the singular cohomology $H^*(Y_{Zar},\mathbb{Q})$. Note that $H^*(Y_{Zar},\mathbb{Q})=H^*(X_K^{an},\mathbb{Q})$, see \cite{bertat}. We show (Theorem \ref{topvan}):\\

{\bf Theorem:} {\it Suppose $Y$ satisfies the monodromy weight conjecture, is of weak Lefschetz type (see section \ref{nundsing}), and that for each $i\ge1$, each component of $Y^i$ is geometrically connected. Suppose $H^i(Y_{Zar},\mathbb{Q})=0$ for all $d>i>0$. Then $N=0$ on $H_{crys}^s(Y/T)_{\mathbb{Q}}$ for all $s\ne d=\dim(Y)$.}\\

We hope that our techniques are useful to further elucidate the constraints which the homotopy type of $X^{an}$ (as Berkovich space) imposes on $N$, and to return to this question in the future. In this connection we mention the work of Le Stum \cite{lestum} dealing with the case of curves. However, in the present paper we work with an abstract strictly semistable log scheme $Y$ over $k$ and only in the final section 7 we consider a lifting $X$ of $Y$ as above.

{\it Notations:} For log algebraic geometry we refer to K. Kato \cite{kalo}. For a (formal) log scheme $(X,{\mathcal N}_X\to{\mathcal O}_X)$ we will often just write $X$ if it is clear from the context to which log structure on $X$ we refer. In this text, all log schemes and morphisms of log schemes have charts for the {\it Zariski} topology. For elements $\{f_j\}_j$ in the structure sheaf of a (formal) scheme we denote by $\mathbb{V}(\{f_j\}_j)$ the closed (formal) subscheme defined by dividing out $\{f_j\}_j$ from the structure sheaf. We let $k$ be a perfect field of characteristic $p>0$ and $W(k)$ its ring of Witt vectors with Frobenius endomorphism $\sigma$. For (sheaves of) $W(k)$-modules $M$ endowed with a $\sigma$-linear endomorphism $F$, and $a\in\mathbb{N}$, we denote by $M(-a)$ the same $W(k)$-module, but now endowed with the endomorphism $p^a.F$. We let $K_0=\quot(W(k))$ and denote by $W(k)\{t\}$ the $p$-adic completion of $W(k)[t]$. We need the formal log scheme $$T=(\spf(W(k)),(\mathbb{N}\longrightarrow W(k), 1\mapsto 0)).$$We denote by $T_1$ its reduction modulo $p$: the log point. For a $p$-adic formal $W(k)$-scheme ${\mathcal F}$ topologically of finite type we denote by ${\mathcal F}_{\mathbb{Q}}$ its generic fibre, as a $K_0$-rigid space. For {\it any} such ${\mathcal F}$ we write $sp$ for the specialization map  ${\mathcal F}_{\mathbb{Q}}\to {\mathcal F}$. For a subscheme $F$ of ${\mathcal F}$ we denote by $]F[_{{\mathcal F}}$ the tube of $F$ in ${\mathcal F}_{\mathbb{Q}}$, i.e. the preimage of $F$ under $sp$; thus $]F[_{{\mathcal F}}$ is an admissible open subspace of ${\mathcal F}_{\mathbb{Q}}$. We use repeatedly and without further comment Kiehl's acyclicity theorem \cite{kiaub} which implies that for a coherent sheaf ${\mathcal M}$ on $]F[_{{\mathcal F}}$ and $i>0$ the push forward sheaves $R^isp_*{\mathcal M}$ on $F$ vanish. 

\section{The tube cohomology of the strata of $Y$}
\label{defy}

\addtocounter{satz}{1}{{{}}\arabic{section}.\arabic{satz}\newcounter{dfsst1}\newcounter{dfsst2}\setcounter{dfsst1}{\value{section}}\setcounter{dfsst2}{\value{satz}} Our basic object of study in this paper is a strictly semistable log scheme $Y$ over $k$. By definition, $Y$ is a fine $T_1$-log scheme $(Y,{\mathcal N}_Y)$ which allows a {\it Zariski} open covering by open subschemes $Y'\subset Y$ with the following property: there exist integers $m\ge 1$ and charts $\mathbb{N}^m\to{\mathcal N}_{Y}(Y')$ for ${\mathcal N}_Y|_{Y'}$ such that\\(i) if on the log scheme $T_1$ we use the chart $\mathbb{N}\to k, 1\mapsto 0$, the diagonal morphism $\mathbb{N}\stackrel{\delta}{\to}\mathbb{N}^m$ is a chart for the structure morphism of log schemes $Y'\to T_1$, and\\(ii) the induced morphism of schemes $$Y'\longrightarrow\spec(k)\times_{\spec(k[t])}\spec(k[t_1,\ldots,t_m])$$ is smooth in the classical sense. If not said otherwise {\it we endow subschemes of} $Y$ {\it with the pull back structure of $T$-log scheme induced by that of} $Y$. By $\{Y_s\}_{s\in R}$ we denote the set of irreducible components of $Y$. {\it We fix an ordering of} $R$. The existence of charts as above for the {\it Zariski}-topology implies that all $Y_j$ are classically smooth. We assume that all connected components of $Y$ are of the same dimension $d$. For $i\in\mathbb{N}$ let $S_i$ be the set of subsets of $R$ with precisely $i$ distinct elements. {\it We identify $S_1$ with $R$}. For $\sigma\in S_i$ let $$Y_{\sigma}^i=\cap_{s\in \sigma}Y_{s}\quad\quad\quad Y^i=\coprod_{\sigma\in S_i}Y_{\sigma}^i.$$(In $Y_{\sigma}^i$ the upper index $i$ is redundant, but it reminds us of the cardinality of $\sigma$)\\

\addtocounter{satz}{1}{{{}} \arabic{section}.\arabic{satz}}\newcounter{adli1}\newcounter{adli2}\setcounter{adli1}{\value{section}}\setcounter{adli2}{\value{satz}} Define the formal log scheme $V=(\spf(W(k)\{t\}),(\mathbb{N}\longrightarrow W(k)\{t\}, 1\mapsto t))$. An {\it admissible lift} of the semistable $k$-log scheme $(Y,{\mathcal N}_Y)$ is a formal $V$-log scheme $({\mathcal Z},{\mathcal N}_{\mathcal Z})$ together with an isomorphism of $T_1$-log schemes $$(Y,{\mathcal N}_Y)\cong({\mathcal Z},{\mathcal N}_{\mathcal Z})\times_{V}T_1$$(where $T_1\to {V}$ is given by $t\mapsto 0$) satisfying the following conditions: On underlying formal schemes ${\mathcal Z}$ is smooth over $\spf(W(k))$, flat over $\spf (W(k)\{t\})$ and its reduction mod $(p)$ is generically smooth over $\spec (k[t])$, the fibre ${\mathcal Y}=\mathbb{V}(t)$ above $t=0$ is a divisor with normal crossings on ${\mathcal Z}$, and ${\mathcal N}_{\mathcal Z}$ is the log structure defined by this divisor. We usually denote an admissible lift by $({\mathcal Z},{\mathcal Y})$. Locally on $Y$, admissible lifts exist. Indeed, by \cite{fkato} 11.3 we locally find embeddings of $Y$ as a normal crossings divisor into smooth $k$-schemes $Z_1$. Assuming $Z_1$ is affine we can lift $Z_1$ to a formally smooth affine formal $W(k)$-scheme ${\mathcal Z}$. Then we lift equations of $Y$ in ${\mathcal O}_{Z_1}$ (which form part of a local system of coordinates on $Z_1$) to equations in ${\mathcal O}_{\mathcal Z}$: these define ${\mathcal Y}$.\\

\addtocounter{satz}{1}{{{}} \arabic{section}.\arabic{satz}}\newcounter{dftuco1}\newcounter{dftuco2}\setcounter{dftuco1}{\value{section}}\setcounter{dftuco2}{\value{satz}} The following construction of diagonal embeddings in the logarithmic context is classical, see for example \cite{hy}, \cite{mokr}. Choose an open covering $Y=\cup_{h\in H}U_h$ of $Y$, together with admissible lifts $({\mathcal Z}_h,{\mathcal Y}_h)$ of the $U_h$ (so $U_h$ is the reduction of ${\mathcal Y}_h$). For a subset $G\subset H$ let $U_G=\cap_{h\in G}U_h$. For $h\in G$ and $s\in R$ we let ${\mathcal Y}_{h,s}$ be the unique $W(k)$-flat irreducible component of ${\mathcal Y}_h$ which lifts $Y_s\cap U_h$; if $Y_s\cap U_h$ is empty we also let ${\mathcal Y}_{h,s}$ be the empty formal scheme. Let ${\mathcal K}''_G$ be the blowing up of $\times_{W}({\mathcal Z}_h)_{h\in G}$ along $\sum_{s\in R}(\times_{W}({\mathcal Y}_{h,s})_{h\in G})$, let ${\mathcal K}'_G$ be the complement of the strict transforms in ${\mathcal K}''_G$ of all ${\mathcal Y}_{h_0,s}\times (\times({\mathcal Z}_h)_{h\in G-\{h_0\}})$ (i.e. all $h_0\in G$, all $s\in R$), and let ${\mathcal Y}'_G$ be the exceptional divisor in ${\mathcal K}'_G$. 

Let $V'_G$ be the blowing up of the $G$-indexed self-product $\times_WV_{h\in G}$ of $V$ along $\times_WT_{h\in G}$; the diagonal embedding $V\to \times_WV_{h\in G}$ lifts to an embedding $V\to V'_G$. There is a natural morphism of formal log schemes ${\mathcal Y}'_G\to V'_G$ and ${\mathcal K}_G={\mathcal K}'_G\times_{V'_G}V$ is a smooth formal $W$-scheme and $V$-log formal scheme with relative normal crossings divisor ${\mathcal Y}_G={\mathcal Y}'_G\times_{V'_G}T$. Given $s\in R$ and any $h\in G$, the closed subscheme ${\mathcal Y}_{G,s}={\mathcal Y}_{h,s}\times_{{\mathcal Y}_h}{\mathcal Y}_G$ of ${\mathcal Y}_G$ is independent of $h$. The non-empty ${\mathcal Y}_{G,s}$ form the set of $W(k)$-flat irreducible components of ${\mathcal Y}_G$. By construction, the diagonal embedding $U_G\to\times_{{W}}({\mathcal Y}_h)_{h\in G}$ lifts canonically to an embedding $U_G\to {\mathcal K}_G$, which in turn factors over the closed formal subscheme ${\mathcal Y}_G$ of ${\mathcal K}_G$. Similarly, if for $i\ge 1$ and $\sigma\in S_i$ we let ${\mathcal Y}^i_{G,\sigma}=\cap_{s\in\sigma}{\mathcal Y}_{G,s}$, then ${Y}^i_{G,\sigma}$ maps to ${\mathcal Y}^i_{G,\sigma}$. 

\begin{lem}\label{dpring} For any $\emptyset \ne G_1\subset G$ the maps$$]U_G[_{{\mathcal Y}_{G}}\longrightarrow]U_{G}[_{{\mathcal Y}_{G_1}},$$$$]U_G\cap Y^i_{\sigma}[_{{\mathcal Y}^i_{G,\sigma}}\longrightarrow]U_{G}\cap Y^i_{\sigma}[_{{\mathcal Y}^i_{G_1,\sigma}}$$are relative open polydisks.\end{lem}

{\sc Proof:} We give a description in local coordinates. Since the statement is local we may assume that there are an $m$ with $1\le m\le i$ and for all $h\in G$ \'{e}tale maps$${\mathcal Z}_h\longrightarrow\spf(W(k)\{t_{h,1},\ldots,t_{h,d}\})$$such that ${\mathcal Y}_h={\mathbb{V}}(t_h)$ with $t_h=\prod_{j=1}^mt_{h,j}$ and such that they induce \'{e}tale maps $$U_G\longrightarrow\spec(k)\times_{\spec(k[t_h])}\spec(k[t_{h,1},\ldots,t_{h,d}])$$which are {\it the same} for all $h$ if (for each fixed $j$) we identify the free variables $t_{h,j}$ for all $h$. We then obtain an \'{e}tale map$${\mathcal K}_G\longrightarrow\spf(\frac{W(k)\{t_{h,1},\ldots,t_{h,d}\}_{h\in G}}{(t_{h}-t_{h'})_{h,h'\in G}}).$$Here the relations $t_{h}=t_{h'}$ are due to the base change $V\to V'_G$ in the definition of ${\mathcal K}_G$. Now since in the definition of ${\mathcal K}'_G$ we removed the strict transforms of all ${\mathcal Y}_{h_0,s}\times (\times({\mathcal Z}_h)_{h\in G-\{h_0\}})$ we may speak of the global sections $t_{h,j}t^{-1}_{h',j}$ in ${\mathcal O}_{{\mathcal K}_G}$ for all $h, h'$. Thus, fixing an element $h_0\in G_1$ we may speak of $v_{h,j}=t_{h,j}t^{-1}_{h_0,j}$ (all $h$). We then have the \'{e}tale map$${\mathcal K}_G\longrightarrow\spf(W(k)\{t_{h_0,1},\ldots,t_{h_0,d},v_{h,2}^{\pm},\ldots,v_{h,d}^{\pm}\}_{h\in G-\{h_0\}})$$(we sold the relations $t_{h}=t_{h'}$ for the price of omitting the terms $v_{h,1}^{\pm}$). Now ${\mathcal Y}_G$ is defined inside ${\mathcal K}_G$ through $t_{h_0}$, and we may arrange the situation in such a way that ${\mathcal Y}^i_{G,\sigma}$ is defined through $t_{h_0,1},\ldots,t_{h_0,i}$. We get the \'{e}tale maps\begin{gather}{\mathcal Y}_G\longrightarrow\spf(W(k)\{t_{h_0,1},\ldots,t_{h_0,d},v_{h,2}^{\pm},\ldots,v_{h,d}^{\pm}\}_{h\in G-\{h_0\}}/(t_{h_0})),\label{polal}\\{\mathcal Y}^i_{G,\sigma}\longrightarrow\spf(W(k)\{t_{h_0,1},\ldots,t_{h_0,d},v_{h,2}^{\pm},\ldots,v_{h,d}^{\pm}\}_{h\in G-\{h_0\}}/(t_{h_0,1},\ldots,t_{h_0,i})).\label{polpa}\end{gather}The closed immersions $U_G\to{\mathcal Y}_G$ and $U_G\cap Y^i_{\sigma}\to{\mathcal Y}^i_{G,\sigma}$ are defined by $p$ and all $v_{h,j}-1$ (all $h\in G-\{h_0\}$, all $2\le j\le d$), hence the map (\ref{polal}) (resp. (\ref{polpa})) induces an isomorphism from $]U_G[_{{\mathcal Y}_{G}}$ (resp. from $]U_G\cap Y^i_{\sigma}[_{{\mathcal Y}^i_{G,\sigma}}$) to the tube over the closed subscheme of the right hand side in (\ref{polal}) (resp. (\ref{polpa})) defined by all $v_{j,h}-1$. To compute this tube in the right hand side we look at the completion along the ideal defined by all $v_{j,h}-1$: it is a formal power series ring over $W(k)\{t_{h_0,1},\ldots,t_{h_0,d}\}/(t_{h_0})$ (resp. over $W(k)\{t_{h_0,i+1},\ldots,t_{h_0,d}\}/(t_{h_0,1},\ldots,t_{h_0,i})$). Repeating all this with the subset $G_1$ of $G$ simply means omitting the free variables $v_{j,h}-1$ for $h\in G-(G_1\cup\{h_0\})$: in both cases the difference is a relative formal powers series ring in the free variables $v_{j,h}-1$ for $h\in G-(G_1\cup\{h_0\})$; but formal power series rings correspond to open polydisks.\\

\addtocounter{satz}{1}{{{}} \arabic{section}.\arabic{satz}}\newcounter{frd1}\newcounter{frd2}\setcounter{frd1}{\value{section}}\setcounter{frd2}{\value{satz}} Endow ${\mathcal K}_G$ with the log structure defined by ${\mathcal Y}_G$, and endow ${\mathcal Y}_G$ with the pull back log structure. Then ${\mathcal K}_G$ (resp. ${\mathcal Y}_G$) is log smooth over $V$ (resp. $T$), and $U_G\to {\mathcal Y}_G$ is an exact closed embedding of (formal) $T$-log schemes. Denote by $\widetilde{\omega}_{{\mathcal K}_G}^{\bullet}$ the logarithmic de Rham complex of ${\mathcal K}_G\to\spf(W(k))$ with the trivial log structure on $\spf(W(k))$. Let $P_{\bullet}\widetilde{\omega}_{{\mathcal K}_G}^{\bullet}$ be the weight filtration on $\widetilde{\omega}_{{\mathcal K}_G}^{\bullet}$: $$P_j\widetilde{\omega}_{{\mathcal K}_G}^{k}=\bi(\widetilde{\omega}_{{\mathcal K}_G}^{j}\otimes\Omega^{k-j}_{{\mathcal K}_G}\longrightarrow\widetilde{\omega}_{{\mathcal K}_G}^{k})$$where $\Omega^{\bullet}_{{\mathcal K}_G}$ is the usual de Rham complex of the morphism of schemes underlying ${\mathcal K}_G\to\spf(W(k))$. For $G_1\subset G_2$ we have natural transition maps ${{\mathcal K}_{G_2}}\to{{\mathcal K}_{G_1}}$, hence a simplicial formal scheme ${\mathcal K}_{\bullet}=\{{\mathcal K}_{G}\}_{G\subset H}$ with sheaf complexes $P_j\widetilde{\omega}_{{\mathcal K}_{\bullet}}^{\bullet}$ on it. Also we have the closed simplicial formal sub schemes ${\mathcal Y}_{\bullet}=\{{\mathcal Y}_{G}\}_{G\subset H}$ and ${\mathcal Y}^i_{\bullet,\sigma}=\{{\mathcal Y}^i_{G,\sigma}\}_{G\subset H}$ for $i\ge 1$ and $\sigma\in S_i$. Denote by ${\mathcal J}_{{\mathcal Y}_{\bullet}}$ (resp. ${\mathcal J}_{{\mathcal Y}^i_{\bullet,\sigma}}$) the ideal of ${\mathcal Y}_{\bullet}$ (resp. of ${\mathcal Y}^i_{\bullet,\sigma}$) in ${\mathcal O}_{{\mathcal K}_{\bullet}}$. Write $\theta=\dlog(t)$. Using the structure sheaf of the simplicial rigid space ${]U_{\bullet}[_{{\mathcal Y}_{\bullet}}}$ we give analytic analogs of the crystalline definitions from \cite{hy}, \cite{mokr}:\begin{align}{\bf C}\widetilde{\omega}^{\bullet}_{{Y}_{\bullet}}&=\frac{\widetilde{\omega}_{{\mathcal K}_{\bullet}}^{\bullet}}{{\mathcal J}_{{\mathcal Y}_{\bullet}}\otimes\widetilde{\omega}_{{\mathcal K}_{\bullet}}^{\bullet}}\otimes_{{\mathcal O}_{{\mathcal Y}_{\bullet}}}sp_*{\mathcal O}_{]U_{\bullet}[_{{\mathcal Y}_{\bullet}}}\notag\\P_j{\bf C}\widetilde{\omega}^{\bullet}_{{Y}_{\bullet}}&=\frac{P_j\widetilde{\omega}_{{\mathcal K}_{\bullet}}^{\bullet}}{{\mathcal J}_{{\mathcal Y}_{\bullet}}\otimes\widetilde{\omega}_{{\mathcal K}_{\bullet}}^{\bullet}}\otimes_{{\mathcal O}_{{\mathcal Y}_{\bullet}}}sp_*{\mathcal O}_{]U_{\bullet}[_{{\mathcal Y}_{\bullet}}}\notag\\{\bf C}{\omega}^{\bullet}_{{Y}_{\bullet}}&=\frac{{\bf C}\widetilde{\omega}^{\bullet}_{{Y}_{\bullet}}}{{\bf C}\widetilde{\omega}^{\bullet-1}_{{Y}_{\bullet}}\wedge\theta}.\notag\end{align}The following definitions which use the structure sheaves of the simplicial rigid spaces ${]U_{\bullet}\cap Y^i_{\sigma}[_{{\mathcal Y}^i_{\bullet,\sigma}}}$ have no analog in \cite{hy}, \cite{mokr}:\begin{align}{\bf C}\widetilde{\omega}^{\bullet}_{{Y}^i_{\bullet,\sigma}}&=\frac{\widetilde{\omega}_{{\mathcal K}_{\bullet}}^{\bullet}}{{\mathcal J}_{{\mathcal Y}^i_{\bullet,\sigma}}\otimes\widetilde{\omega}_{{\mathcal K}_{\bullet}}^{\bullet}}\otimes_{{\mathcal O}_{{\mathcal Y}^i_{\bullet,\sigma}}}sp_*{\mathcal O}_{]U_{\bullet}\cap Y^i_{\sigma}[_{{\mathcal Y}^i_{\bullet,\sigma}}}\notag\\P_j{\bf C}\widetilde{\omega}^{\bullet}_{{Y}^i_{\bullet,\sigma}}&=\frac{P_j\widetilde{\omega}_{{\mathcal K}_{\bullet}}^{\bullet}+{\mathcal J}_{{\mathcal Y}^i_{\bullet,\sigma}}\otimes\widetilde{\omega}_{{\mathcal K}_{\bullet}}^{\bullet}}{{\mathcal J}_{{\mathcal Y}^i_{\bullet,\sigma}}\otimes\widetilde{\omega}_{{\mathcal K}_{\bullet}}^{\bullet}}\otimes_{{\mathcal O}_{{\mathcal Y}^i_{\bullet,\sigma}}}sp_*{\mathcal O}_{]U_{\bullet}\cap Y^i_{\sigma}[_{{\mathcal Y}^i_{\bullet,\sigma}}}\notag\\{}&=\frac{P_j\widetilde{\omega}_{{\mathcal K}_{\bullet}}^{\bullet}}{{\mathcal J}_{{\mathcal Y}^i_{\bullet,\sigma}}\otimes\widetilde{\omega}_{{\mathcal K}_{\bullet}}^{\bullet}\cap P_j\widetilde{\omega}_{{\mathcal K}_{\bullet}}^{\bullet}}\otimes_{{\mathcal O}_{{\mathcal Y}^i_{\bullet,\sigma}}}sp_*{\mathcal O}_{]U_{\bullet}\cap Y^i_{\sigma}[_{{\mathcal Y}^i_{\bullet,\sigma}}}\notag\\{\bf C}{\omega}^{\bullet}_{{Y}^i_{\bullet,\sigma}}&=\frac{{\bf C}\widetilde{\omega}^{\bullet}_{{Y}^i_{\bullet,\sigma}}}{{\bf C}\widetilde{\omega}^{\bullet-1}_{{Y}^i_{\bullet,\sigma}}\wedge\theta}\notag\\{\bf C}\Omega^{\bullet}_{Y_{\bullet,\sigma}^i}&=P_0{\bf C}\widetilde{\omega}^{\bullet}_{Y_{\bullet,\sigma}^i}.\notag\end{align}We drop the $\sigma$ in these notations when we sum over all $\sigma\in S_i$:$${\bf C}\widetilde{\omega}^{\bullet}_{Y_{\bullet}^i}=\bigoplus_{\sigma\in S_i}{\bf C}\widetilde{\omega}^{\bullet}_{Y_{\bullet,\sigma}^i}\quad\quad\quad P_j{\bf C}\widetilde{\omega}^{\bullet}_{Y_{\bullet}^i}=\bigoplus_{\sigma\in S_i}P_j{\bf C}\widetilde{\omega}^{\bullet}_{Y_{\bullet,\sigma}^i},$$$${\bf C}{\omega}^{\bullet}_{Y_{\bullet}^i}=\bigoplus_{\sigma\in S_i}{\bf C}{\omega}^{\bullet}_{Y_{\bullet,\sigma}^i}\quad\quad\quad{\bf C}\Omega^{\bullet}_{Y_{\bullet}^i}=\bigoplus_{\sigma\in S_i}{\bf C}\Omega^{\bullet}_{Y_{\bullet,\sigma}^i}.$$ We view all these sheaf complexes as living on $U_{\bullet}=\{U_G\}_G$. In our notation we will frequently drop the subscript bullet below $Y$ (which holds the place for the varying $G$), thus we understand$${\bf C}\widetilde{\omega}^{\bullet}_{{Y}}={\bf C}\widetilde{\omega}^{\bullet}_{{Y}_{\bullet}},\quad\quad\quad{\bf C}{\omega}^{\bullet}_{{Y}}={\bf C}{\omega}^{\bullet}_{{Y}_{\bullet}},$$$${\bf C}\widetilde{\omega}^{\bullet}_{{Y}_{\sigma}^i}={\bf C}\widetilde{\omega}^{\bullet}_{{Y}^i_{\bullet,\sigma}},\quad\quad\quad{\bf C}{\omega}^{\bullet}_{{Y}_{\sigma}^i}={\bf C}\widetilde{\omega}^{\bullet}_{{Y}^i_{\bullet,\sigma}},$$$${\bf C}\widetilde{\omega}^{\bullet}_{{Y}^i}={\bf C}\widetilde{\omega}^{\bullet}_{{Y}^i_{\bullet}},\quad\quad\quad{\bf C}{\omega}^{\bullet}_{{Y}^i}={\bf C}{\omega}^{\bullet}_{{Y}^i_{\bullet}}.$$Moreover, to be consistent with the introduction we keep the names of complexes on $U_{\bullet}$ also for their derived push forward on $Y$ (via the morphism of simplicial schemes $U_{\bullet}\to Y$). Note that ${\bf C}\Omega^{\bullet}_{Y_{\sigma}^i}$ computes the non logarithmic convergent (or equivalently: crystalline) cohomology of the classically smooth $k$-scheme $Y_{\sigma}^i$.\\
Recall from \cite{hyoka} 3.1 that on the formal log scheme $T$ we have a Frobenius action: the unique endomorphism which equals $\sigma$ on $W(k)$ and multiplication by $p$ on the standard chart $\mathbb{N}\stackrel{0}{\to}W(k)$. We may assume that for each $h\in H$ there is an endomorphism of $({\mathcal Z}_h,{\mathcal Y}_h)$ which lifts the Frobenius endomorphism (i.e. the $p$-power map on the structure sheaf) of its reduction modulo $(p)$ and which sends equations for the divisors ${\mathcal Y}_{h,s}$ on ${\mathcal Z}_h$ to their $p$-th power. Then we also get such Frobenius endomorphisms $F$ on ${\mathcal K}_{\bullet}$, ${\mathcal Y}_{\bullet}$ and ${\mathcal Y}^i_{\bullet,\sigma}$ and on the various logarithmic de Rham complexes defined above. By abuse of notation we will frequently write $\mathbb{R}\Gamma(Y,K^{\bullet})$ instead of $\mathbb{R}\Gamma(U_{\bullet},K^{\bullet})$ for sheaf complexes $K^{\bullet}$ on $U_{\bullet}$.

\begin{pro}\label{indeplif} The objects $\mathbb{R}\Gamma(Y,{\bf C}{\omega}^{\bullet}_{{Y}_{\bullet}})$, $\mathbb{R}\Gamma(Y,P_j{\bf C}\widetilde{\omega}^{\bullet}_{{Y}_{\bullet}})$, $\mathbb{R}\Gamma(Y,{\bf C}{\omega}^{\bullet}_{{Y}_{\bullet,\sigma}^i})$ and $\mathbb{R}\Gamma(Y,P_j{\bf C}\widetilde{\omega}^{\bullet}_{{Y}_{\bullet,\sigma}^i})$ are independent of the chosen system $\{({\mathcal Z}_h,{\mathcal Y}_h)\}_{h\in H}$.
\end{pro}

{\sc Proof:} Given another system $\{({\mathcal Z}'_{h'},{\mathcal Y}'_{h'})\}_{h'\in H'}$ one performs the constructions from \arabic{dftuco1}.\arabic{dftuco2} also for $\{({\mathcal Z}'_{h'},{\mathcal Y}'_{h'})\}_{h'\in H'}$ and for the union of the systems $\{({\mathcal Z}_h,{\mathcal Y}_h)\}_{h\in H}$ and $\{({\mathcal Z}'_{h'},{\mathcal Y}'_{h'})\}_{h'\in H'}$. We get canonical maps from the cohomology object formed with respect to this union to those formed with respect to $\{({\mathcal Z}_h,{\mathcal Y}_h)\}_{h\in H}$ and $\{({\mathcal Z}'_{h'},{\mathcal Y}'_{h'})\}_{h'\in H'}$. That these are isomorphisms is a local statement and follows from \ref{dpring} and the Poincar\'{e} lemma for relative open polydisks.\\

\addtocounter{satz}{1}{{{}} \arabic{section}.\arabic{satz}} The logarithmic de Rham complex $\omega^{\bullet}_{{\mathcal Y}_{\bullet}/T}$ of ${\mathcal Y}_{\bullet}\to T$ is $$\omega^{\bullet}_{{\mathcal Y}_{\bullet}/T}=\frac{\widetilde\omega^{\bullet}_{{\mathcal Y}_{\bullet}/T}}{\widetilde\omega^{\bullet-1}_{{\mathcal Y}_{\bullet}/T}\wedge\theta}\quad\quad\mbox{with}\quad\quad\widetilde\omega^{\bullet}_{{\mathcal Y}_{\bullet}/T}=\widetilde{\omega}_{{\mathcal K}_{\bullet}}^{\bullet}\otimes_{{\mathcal O}_{{\mathcal K}_{\bullet}}}{\mathcal O}_{{\mathcal Y}_{\bullet}}.$$The logarithmic convergent cohomology in our context (see Ogus \cite{ogcon} and Shiho \cite{shiho} for more general definitions) is given by the objects$$\mathbb{R}\Gamma_{conv}(Y/T)=\mathbb{R}\Gamma(Y,\omega^{\bullet}_{{\mathcal Y}_{\bullet}/T}\otimes_{{\mathcal O}_{{\mathcal Y}_{\bullet}}}sp_*{\mathcal O}_{{]U_{\bullet}[_{{\mathcal Y}_{\bullet}}}}),$$$$\mathbb{R}\Gamma_{conv}(Y^i_{\sigma}/T)=\mathbb{R}\Gamma(Y,\omega^{\bullet}_{{\mathcal Y}_{\bullet}/T}\otimes_{{\mathcal O}_{{\mathcal Y}_{\bullet}}}sp_*{\mathcal O}_{{]U_{\bullet}\cap Y^i_{\sigma}[_{{\mathcal Y}_{\bullet}}}}).$$

\begin{pro}\label{cocri} (i) We have a canonical isomorphism$$\mathbb{R}\Gamma_{conv}(Y/T)\cong \mathbb{R}\Gamma_{crys}(Y/T)_{\mathbb Q}.$$(ii) There are canonical isomorphisms$$\mathbb{R}\Gamma_{conv}(Y^i_{\sigma}/T)\cong\mathbb{R}\Gamma(Y,{\bf C}{\omega}^{\bullet}_{{Y}_{\bullet,\sigma}^i}).$$\end{pro}

{\sc Proof:} (i) is due to a general comparison isomorphism between log convergent and log crystalline cohomology of a log smooth morphism, see e.g. \cite{shiho}. In (ii) the canonical map is induced by the inclusion of simplicial rigid spaces ${]U_{\bullet}\cap Y^i_{\sigma}[_{{\mathcal Y}^i_{\bullet,\sigma}}}\to{]U_{\bullet}\cap Y^i_{\sigma}[_{{\mathcal Y}_{\bullet}}}$. That it is an isomorphism can be checked locally, so we may assume that there exists an affine admissible lift $({\mathcal Z},{\mathcal Y})$ of $Y$. For $s\in R$ let ${\mathcal Y}_{s}$ be the $W(k)$-flat irreducible component of ${\mathcal Y}$ lifting $Y_{s}$. For $\tau\subset R$ let ${\mathcal Y}_{\tau}=\cap_{s\in\tau}{\mathcal Y}_{s}$ if $\emptyset\ne\tau$ and ${\mathcal Y}_{\emptyset}={\mathcal Y}$. Then let $${\mathcal O}_{\tau}=\Gamma(]Y_{\tau\cup\sigma}^{|\tau\cup\sigma|}[_{{\mathcal Y}_{\tau}},{\mathcal O}_{{\mathcal Y}_{\tau}}).$$Let $\omega^{\bullet}_{\mathcal Y}$ be the logarithmic de Rham complex of ${\mathcal Y}/T$ and write $L^{\bullet}=\Gamma(Y,\omega^{\bullet}_{\mathcal Y}\otimes\mathbb{Q})$. Then$$\mathbb{R}^*\Gamma_{conv}(Y^i_{\sigma}/T)=h^*(L^{\bullet}\otimes{\mathcal O}_{\emptyset})$$$$\mathbb{R}\Gamma(Y,{\bf C}{\omega}^{\bullet}_{{Y}_{\bullet,\sigma}^i})=h^*(L^{\bullet}\otimes{\mathcal O}_{\sigma}).$$Thus we need to show that $L^{\bullet}\otimes{\mathcal O}_{\emptyset}\to L^{\bullet}\otimes{\mathcal O}_{\sigma}$ is a quasiisomorphism. The following purely formal reduction to Lemma \ref{schritt} below is literally the same as in \cite{hkstrat} Proposition 4.2 or \cite{colo} Theorem 3.14. For subsets $\mu\subset R$ we may form the closed formal subscheme ${\mathcal Y}^{\mu}=\cup_{s\in\mu}{\mathcal Y}_{s}$ of ${\mathcal Y}$ (not to be confused with our previous notation ${\mathcal Y}^j$ for $j\in\mathbb{N}$) and write $${\mathcal O}^{\mu}=\Gamma(]Y_{\sigma}^{i}\cap{\mathcal Y}^{\mu}[_{{\mathcal Y}^{\mu}},{\mathcal O}_{{\mathcal Y}^{\mu}}).$$Similarly, for two subsets $\mu_1, \mu_2$ of $R$ we write $${\mathcal O}^{\mu_1,\mu_2}=\Gamma(]Y_{\sigma}^{i}\cap{\mathcal Y}^{\mu_1}\cap{\mathcal Y}^{\mu_2}[_{{\mathcal Y}^{\mu_1}\cap{\mathcal Y}^{\mu_2}},{\mathcal O}_{{\mathcal Y}^{\mu_1}\cap{\mathcal Y}^{\mu_2}})$$and$${\mathcal O}^{\mu_1}_{\mu_2}=\Gamma(]Y_{\sigma\cup\mu_2}^{|\sigma\cup\mu_2|}\cap{\mathcal Y}^{\mu_1}[_{{\mathcal Y}^{\mu_1}\cap{\mathcal Y}_{\mu_2}},{\mathcal O}_{{\mathcal Y}^{\mu_1}\cap{\mathcal Y}_{\mu_2}}).$$We will show that in$$L^{\bullet}\otimes{\mathcal O}_{\emptyset}\stackrel{\alpha}{\longrightarrow} L^{\bullet}\otimes{\mathcal O}^{\sigma}\stackrel{\beta}{\longrightarrow}L^{\bullet}\otimes{\mathcal O}_{\sigma}$$both $\alpha$ and $\beta$ are quasiisomorphisms. The exact sequences$$0\longrightarrow {\mathcal O}_{\emptyset}\longrightarrow{\mathcal O}^{\sigma}\oplus{\mathcal O}^{R-\sigma}\longrightarrow{\mathcal O}^{\sigma, R-\sigma}\longrightarrow0$$$$0\longrightarrow{\mathcal O}^{\sigma}\longrightarrow{\mathcal O}^{\sigma}\oplus{\mathcal O}^{\sigma, R-\sigma}\longrightarrow{\mathcal O}^{\sigma, R-\sigma}\longrightarrow0$$ show that, to prove that $\alpha$ is a quasiisomorphism, it is enough to prove that $L^{\bullet}\otimes{\mathcal O}^{R-\sigma}\to L^{\bullet}\otimes{\mathcal O}^{\sigma, R-\sigma}$ is a quasiisomorphism. To see this, it is enough to show that both $L^{\bullet}\otimes{\mathcal O}^{R-\sigma}\stackrel{\gamma}{\to}L^{\bullet}\otimes{\mathcal O}^{R-\sigma}_{\sigma}$ and $L^{\bullet}\otimes{\mathcal O}^{\sigma, R-\sigma}\stackrel{\delta}{\to}L^{\bullet}\otimes{\mathcal O}^{R-\sigma}_{\sigma}$ are quasiisomorphisms. Consider the exact sequence\begin{gather}0\longrightarrow {\mathcal O}^{R-\sigma}\longrightarrow\bigoplus_{s\in R-\sigma}{\mathcal O}_{s}\longrightarrow\bigoplus_{{\rho\subset R-\sigma}\atop{|\rho|=2}}{\mathcal O}_{\rho}\longrightarrow\ldots\longrightarrow{\mathcal O}_{R-\sigma} \longrightarrow 0\tag{$*$}\end{gather}Comparison of the exact sequences $(*)\otimes L^{\bullet}$ and $(*)\otimes {\mathcal O}^{R-\sigma}_{\sigma}\otimes L^{\bullet}$ shows that to prove that $\gamma$ is a quasiisomorphism, it is enough to show this for $L^{\bullet}\otimes{\mathcal O}_{\rho}\to L^{\bullet}\otimes{\mathcal O}_{\rho\cup\sigma}$ for all $\emptyset\ne \rho\subset R-\sigma$; but this is Lemma \ref{schritt}. Comparison of $(*)\otimes L^{\bullet}\otimes{\mathcal O}^{\sigma, R-\sigma}$ and $(*)\otimes L^{\bullet}\otimes{\mathcal O}^{R-\sigma}_{\sigma}$ shows that to prove that $\delta$ is a quasiisomorphism, it is enough to show this for $L^{\bullet}\otimes{\mathcal O}^{\sigma}_{\rho}\stackrel{\epsilon_{\rho}}{\to}L^{\bullet}\otimes {\mathcal O}_{\rho\cup\sigma}$ for all $\emptyset\ne \rho\subset R-\sigma$. Consider the exact sequence \begin{gather}0\longrightarrow {\mathcal O}^{\sigma}\longrightarrow\bigoplus_{s\in \sigma}{\mathcal O}_{s}\longrightarrow\bigoplus_{{\gamma\subset \sigma}\atop{|\gamma|=2}}{\mathcal O}_{\gamma}\longrightarrow\ldots\longrightarrow{\mathcal O}_{\sigma} \longrightarrow 0\tag{$**$}\end{gather} The exact sequence $(**)\otimes {\mathcal O}^{\sigma}_{\rho}\otimes L^{\bullet}$ shows that to prove that $\epsilon_{\rho}$ is a quasiisomorphism, it is enough to show this for $L^{\bullet}\otimes {\mathcal O}_{\rho\cup\gamma}\to L^{\bullet}\otimes {\mathcal O}_{\rho\cup\sigma}$ for all $\emptyset\ne \gamma\subset \sigma$; but this is Lemma \ref{schritt}. The exact sequence $(**)\otimes L^{\bullet}$ shows that to prove that $\beta$ is a quasiisomorphism, it is enough to show this for $L^{\bullet}\otimes {\mathcal O}_{\gamma}\to L^{\bullet}\otimes {\mathcal O}_{\sigma}$
for all $\emptyset\ne \gamma\subset\sigma$; but this is Lemma \ref{schritt}.

\begin{lem}\label{schritt} For any inclusion $\emptyset\ne\tau_1\subset\tau_2\subset R$ with $\tau_2-\tau_1\subset \sigma$ the projection$$\mu:L^{\bullet}\otimes{\mathcal O}_{\tau_1}\longrightarrow L^{\bullet}\otimes{\mathcal O}_{\tau_2}$$is a quasiisomorphism.\end{lem}

{\sc Proof:} Also this is as in \cite{hkstrat} Proposition 4.2. By induction we may suppose $\tau_2=\tau_1\cup\{s_0\}$ for some $s_0\in\sigma$ with $s_0\notin \tau_2$. We may assume that there is a smooth morphism $${\mathcal Y}\longrightarrow\spf(W(k))\times_{\spf(W(k)\{t\})}\spf(W(k)\{t_1,\ldots,t_m\})$$lifting the situation described in \arabic{dfsst1}.\arabic{dfsst2}, and furthermore that ${\mathcal Y}_{\tau_2}$ is the closed formal subscheme of ${\mathcal Y}_{\tau_1}$ defined by $U:=t_1$ (both are closed formal subschemes of ${\mathcal Y}$). Then$${\mathcal O}_{\tau_1}=\{\sum_{m\ge0}a_mU^m;\,a_m\in{\mathcal O}_{\tau_2}, \ord_p(a_m)\longrightarrow\infty\}.$$Localizing further we may assume that there exist $s_1,\ldots,s_n\in{\mathcal O}_{{\mathcal Y}}({\mathcal Y})$ (with $n+m-1=\dim(Y)$) such that $\{ds_1,\ldots,ds_n,\dlog(t_1),\ldots,\dlog(t_{m-1})\}$ is an ${\mathcal O}_{\mathcal Y}({\mathcal Y})$-basis of $L^1$. The images of these elements (denoted by the same names) in $L^1\otimes{\mathcal O}_{\tau_i}$ then form an ${\mathcal O}_{\tau_i}$-basis of $L^1\otimes{\mathcal O}_{\tau_i}$ ($i=1,2$). Let $L_c^1$, resp. $<\dlog(U)>$, be the ${\mathcal O}_{\tau_2}$-submodule of $L^1\otimes{\mathcal O}_{\tau_2}$ generated by the set $\{ds_1,\ldots,ds_n,\dlog(t_2),\ldots,\dlog(t_{m-1})\}$, resp. the single element $\dlog(U)=\dlog(t_1)$. Let $L_c^{\bullet}$, resp. $<\dlog(U)>^{\bullet}$ be the sub-${\mathcal O}_{\tau_2}$-algebra of $L^{\bullet}\otimes {\mathcal O}_{\tau_2}$ generated by $L_c^1$, resp. by $<\dlog(U)>$. These are in fact sub{\it complexes} and we have the decomposition$$L^{\bullet}\otimes {\mathcal O}_{\tau_2}=L_c^{\bullet}\otimes<\dlog(U)>^{\bullet}.$$By the above description of ${\mathcal O}_{\tau_1}$ it follows that the map $\mu$ in question has a natural "zero"-section $\nu$ and it suffices to show that $\nu$ induces surjective maps in cohomology. Let $\omega\in L^k\otimes {\mathcal O}_{\tau_1}$. It can be written as$$\omega=\sum_{m\ge 0}a_{m}U^m\dlog(U)+\sum_{m\ge0}b_{m}U^{m}$$with $a_{m}\in L_c^{k-1}$ and $b_{m}\in L_c^k$. Subtracting $d(\sum_{m>0}m^{-1}a_{m}U^{m})$ and renaming the coefficients we may write $\omega$ modulo exact forms as$$\omega=a_0\dlog(U)+\sum_{m\ge0}b_{m}U^{m}.$$If $d\omega=0$ we get $\omega=a_0\dlog(U)+b_0$ which lies in the image of $\nu$.\\

\addtocounter{satz}{1}{{{}} \arabic{section}.\arabic{satz}}\newcounter{idlo1}\newcounter{idlo2}\setcounter{idlo1}{\value{section}}\setcounter{idlo2}{\value{satz}} Remarks. (1) An idealized log scheme is a log scheme together with an ideal in its log structure which maps to the zero element of the structure sheaf, see Ogus \cite{ogid}. There is a notion of ideally log smooth morphisms between idealized log schemes, defined as usual by a lifting condition over exact closed nilimmersions. For $\sigma\in S_i$ let ${\mathcal F}_{i,\sigma}\subset{\mathcal N}_Y$ be the preimage of $\ke({\mathcal O}_Y\to{\mathcal O}_{Y_{\sigma}^i})$ in the log structure ${\mathcal N}_Y$ of $Y$. On $Y_{\sigma}^i$ the pair $({\mathcal N}_Y,{\mathcal F}_{i,{\sigma}})$ induces the structure of an idealized log scheme. If we view $T_1$ as an idealized log scheme by simply taking the zero ideal in its log structure, the morphism $Y_{\sigma}^i\to T_1$ becomes ideally log smooth. One may define the idealized log crystalline site of $Y_{\sigma}^i/T$ and will find that the cohomology (tensored with ${\mathbb Q}$) of its structure sheaf is isomorphic to our $\mathbb{R}\Gamma(Y,{\bf C}{\omega}^{\bullet}_{{Y}_{\bullet,\sigma}^i})$. There are also Cartier isomorphisms valid in this context, one may therefore define corresponding de Rham-Witt complexes and everything we develop here could be done with them as well.\\(2) We could similarly define the idealized log convergent cohomology of $Y_{\sigma}^i/T$: it is isomorphic to the non-idealized version, as follows from \ref{cocri}. (We do not claim coincidence of the idealized and the non-idealized log {\it crystalline} cohomology of $Y_{\sigma}^i/T$; we suspect that the latter is not a useful theory.)\\

\begin{pro}\label{logresi} For $k\ge r\ge 1$ let $\tau\in S_k$, $\rho\in S_r$ such that $\rho\subset\tau$. Taking the residue along ${\mathcal Y}_{\bullet,\rho}^r$ defines a map $Res_{\rho}:\widetilde{\omega}_{{\mathcal K}_{\bullet}}^{q}\to\widetilde{\omega}_{{\mathcal K}_{\bullet}}^{q-r}\otimes{\mathcal O}_{{\mathcal Y}^r_{\bullet,\rho}}$. For $j\ge r$ it extends to a map$$Res_{\rho}:P_j{\bf C}\widetilde{\omega}^{q}_{{Y}^k_{\bullet,\tau}}\longrightarrow P_{j-r}{\bf C}\widetilde{\omega}^{q-r}_{{Y}^k_{\bullet,\tau}}.$$
\end{pro}

{\sc Proof:} Taking residues involves the choice of local coordinates, the problem is to show the independence of this choice. We must work on each ${\mathcal K}_{G}$ (for $G\subset H$) separately. Let us write $\rho=\{s_1,\ldots,s_r\}$ with $s_1<\ldots<s_r$ in our fixed ordering of $R$. We assume that ${\mathcal Y}^r_{G,\rho}\ne\emptyset$. Since we work locally we may suppose that there exist $t_1,\ldots,t_r\in{\mathcal O}_{{\mathcal K}_{G}}$ such that ${\mathcal Y}_{G,s_i}=\mathbb{V}(t_{i})$ for $1\le i\le r$. Let$$\dlog(t_{\rho})=\dlog(t_{r})\wedge\ldots\wedge\dlog(t_{1}).$$For $\mu\subset R$ denote by $\widetilde{\omega}_{{\mathcal K}_{G},\mu}^{\bullet}$ the logarithmic differential module on ${\mathcal K}_{G}$ with logarithmic poles along $\cup_{s\in \mu}{\mathcal Y}_{G,s}$. Let $${\mathcal P}^q=\bi[\widetilde{\omega}_{{\mathcal K}_{G}}^{r-1}\otimes\widetilde{\omega}_{{\mathcal K}_{G},R-\rho}^{q+1-r}\longrightarrow\widetilde{\omega}_{{\mathcal K}_{G}}^{q}].$$Now let $\omega\in\widetilde{\omega}_{{\mathcal K}_{G}}^{q}$. It can be written as $$\omega=\omega_0+\eta\wedge\dlog(t_{\rho})$$with $\omega_0\in{\mathcal P}^q$ and $\eta\in \widetilde{\omega}_{{\mathcal K}_{G},R-\rho}^{q-r}$. We set$$Res_{\rho}(\omega)=\eta\otimes 1\in\widetilde{\omega}_{{\mathcal K}_{G}}^{q-r}\otimes{\mathcal O}_{{\mathcal Y}^r_{G,\rho}}.$$To see that this is well defined consider first another sum decomposition $\omega=\omega'_0+\eta'\wedge\dlog(t_{\rho})$ with $\omega'_0\in{\mathcal P}^q$ and $\eta'\in \widetilde{\omega}_{{\mathcal K}_{G},R-\rho}^{q-r}$. Then $$(\eta-\eta')\wedge\dlog(t_{\rho})=\omega'_0-\omega_0\in{\mathcal P}^q$$which implies that at least one of $t_1,\ldots,t_r$ divides $\eta-\eta'$, hence $\eta\otimes 1=\eta'\otimes 1$ in $\widetilde{\omega}_{{\mathcal K}_{G}}^{q-r}\otimes{\mathcal O}_{{\mathcal Y}^r_{G,\rho}}$.

To see independence of the chosen system $t_1,\ldots,t_r$ consider another system $t'_1,\ldots,t'_r\in{\mathcal O}_{{\mathcal K}_{G}}$ such that ${\mathcal Y}_{G,s_i}=\mathbb{V}(t'_{i})$ for $1\le i\le r$. Any such system $t'_1,\ldots,t'_r$ arises from $t_1,\ldots,t_r$  by finitely many operations of the following type: multiply $t_i$ for a single $1\le i\le r$ by a unit in ${\mathcal O}_{{\mathcal K}_{G}}$. Thus we may assume $t'_{i_0}=\epsilon.t_{i_0}$ for some $1\le i_0\le r$ and $\epsilon\in{{\mathcal O}_{{\mathcal K}_{G}}^{\times}}$, and $t_i'=t_i$ for all $i\ne i_0$. We want to show $Res'_{\rho}(\omega)=Res_{\rho}(\omega)$ for the residue map $Res'_{\rho}$ defined with respect to $\{t'_i\}_{1\le i\le r}$. Write $$\omega=\omega_0+\eta\wedge\dlog(t'_{\rho})+\eta\wedge(\dlog(t_{\rho})-\dlog(t'_{\rho})).$$Clearly $Res'_{\rho}(\omega_0+\eta\wedge\dlog(t'_{\rho}))=\eta$. Therefore we need to show $Res'_{\rho}(\mu)=0$ for $\mu=\eta\wedge(\dlog(t_{\rho})-\dlog(t'_{\rho}))$. Now$$-\mu=\eta\wedge\dlog(t_{r})\wedge\ldots\wedge\dlog(t_{{i_0}+1})\wedge\dlog(\epsilon)\wedge\dlog(t_{i_0-1})\wedge\ldots\wedge\dlog(t_{1})$$and $\dlog(\epsilon)\in P_0\widetilde{\omega}_{{\mathcal K}_{G}}^{1}$, hence $\mu\in {\mathcal P}^q$ and $Res'_{\rho}(\mu)=0$.

\begin{pro}\label{corand} For $i,j\ge1$, $\sigma\in S_i$ there is a canonical isomorphism$$\Gr _j{\bf C}\widetilde{\omega}^{\bullet}_{{Y}^i_{\bullet,\sigma}}\cong\bigoplus_{\tau\in S_j}{\bf C}{\Omega}^{\bullet}_{{Y}^{|\sigma\cup\tau|}_{\bullet,\sigma\cup\tau}}[-j](-j).$$
\end{pro}

{\sc Proof:} (Given $\sigma\in S_i$ and $\tau\in S_j$ we may form the union $\sigma\cup\tau\subset R$, an element of $S_{|\sigma\cup\tau|}$ with $|\sigma\cup\tau|\le |\sigma|+|\tau|=i+j$.) Recall that ${\mathcal Y}_{\bullet}$ is a relative normal crossings divisor in the smooth formal simplicial $W(k)$-scheme ${\mathcal K}_{\bullet}$ and that $\{{\mathcal Y}_{{\bullet},\tau}^j\}_{\tau\in S_j}$ is the set of its $j$-codimensional intersection strata. In such a situation it is a classical fact that taking Poincar\'{e} residues (our \ref{logresi} in the extreme case $j=r$) induces an isomorphism$$ \Gr_j\widetilde{\omega}^{\bullet}_{{\mathcal K}_{{\bullet}}}\cong\bigoplus_{\tau\in S_j}{\Omega}^{\bullet}_{{\mathcal Y}^{j}_{{\bullet},\tau}}[-j](-j)$$where $({\Omega}^{\bullet}_{{\mathcal Y}^{j}_{{\bullet},\tau}},d)$ denotes the classical de Rham complex on the smooth formal simplicial $W(k)$-scheme ${\mathcal Y}^{j}_{{\bullet},\tau}$. We claim that this isomorphism restricts to an isomorphism between the ${\mathcal O}_{{\mathcal Y}_{\bullet}}$-submodule generated by ${\mathcal J}_{{\mathcal Y}_{\bullet,\sigma}^i}\otimes\widetilde{\omega}^{\bullet}_{{\mathcal K}_{\bullet}}\cap P_j\widetilde{\omega}^{\bullet}_{{\mathcal K}_{\bullet}}$ and the ${\mathcal O}_{{\mathcal Y}_{\bullet}}$-submodule generated by $({\Omega}^{\bullet-1}_{{\mathcal Y}^{j}_{{\bullet},\tau}}\wedge d({\mathcal J}_{{\mathcal Y}_{\bullet,\sigma}^i})+{\mathcal J}_{{\mathcal Y}_{\bullet,\sigma}^i}.{\Omega}^{\bullet}_{{\mathcal Y}^{j}_{{\bullet},\tau}})[-j]$. Indeed, we may work locally around ${\mathcal Y}_{\bullet,\sigma\cup\tau}^{|\sigma\cup\tau|}$ for ${\tau\in S_j}$ and assume ${\mathcal Y}_{\bullet,\sigma}^i={\mathbb V}(t_s)_{s\in\sigma}$ and ${\mathcal Y}_{\bullet,\tau}^j={\mathbb V}(t_s)_{s\in\tau}$ for suitable $\{t_s\}_{s\in \sigma\cup\tau}\in{\mathcal O}_{{\mathcal K}_{\bullet}}$. Then the ${\mathcal O}_{{\mathcal Y}_{\bullet}}$-submodule of the above left hand side generated by ${\mathcal J}_{{\mathcal Y}_{\bullet,\sigma}^i}\otimes\widetilde{\omega}^{\bullet}_{{\mathcal K}_{\bullet}}\cap P_j\widetilde{\omega}^{\bullet}_{{\mathcal K}_{\bullet}}$ is in fact generated by elements of the form$$\omega=t_{s_0}\omega'\wedge\bigwedge_{s\in\tau}\dlog(t_s)$$for some $s_0\in\sigma$ and $\omega'\in\widetilde{\omega}^{\bullet-j}_{{\mathcal K}_{{\bullet}}}$ such that $t_{s_0}\omega'\in P_0\widetilde{\omega}^{\bullet-j}_{{\mathcal K}_{{\bullet}}}$. If $\dlog(t_{s_0})$ divides $\omega'$ (in the graded algebra $\widetilde{\omega}^{\bullet}_{{\mathcal K}_{{\bullet}}}$) then $dt_{s_0}$ divides $t_{s_0}\omega'$ and $Res_{\tau}(\omega)\in {\Omega}^{\bullet-1}_{{\mathcal Y}^{j}_{{\bullet},\tau}}\wedge d({\mathcal J}_{{\mathcal Y}_{\bullet,\sigma}^i})$. If $\dlog(t_{s_0})$ does not divide $\omega'$ then $\omega'\in P_0\widetilde{\omega}^{\bullet-j}_{{\mathcal K}_{{\bullet}}}$ and $Res_{\tau}(\omega)\in {\mathcal J}_{{\mathcal Y}_{\bullet,\sigma}^i}.{\Omega}^{\bullet}_{{\mathcal Y}^{j}_{{\bullet},\tau}}$. Conversely, these local considerations show how to fabricate preimages of elements of $({\Omega}^{\bullet-1}_{{\mathcal Y}^{j}_{{\bullet},\tau}}\wedge d({\mathcal J}_{{\mathcal Y}_{\bullet,\sigma}^i})+{\mathcal J}_{{\mathcal Y}_{\bullet,\sigma}^i}.{\Omega}^{\bullet}_{{\mathcal Y}^{j}_{{\bullet},\tau}})[-j]$ under the above isomorphism. The claim is established.

Dividing out these submodules we get the isomorphism$$\frac{P_j\widetilde{\omega}^{\bullet}_{{\mathcal K}_{{\bullet}}}}{P_{j-1}\widetilde{\omega}^{\bullet}_{{\mathcal K}_{{\bullet}}}+({\mathcal J}_{{\mathcal Y}_{\bullet,\sigma}^i}\otimes\widetilde{\omega}^{\bullet}_{{\mathcal K}_{\bullet}}\cap P_j\widetilde{\omega}^{\bullet}_{{\mathcal K}_{\bullet}})}\cong\bigoplus_{\tau\in S_j}{\Omega}^{\bullet}_{{\mathcal Y}^{|\sigma\cap \tau|}_{{\bullet},\sigma\cap\tau}}[-j](-j).$$Tensoring over ${\mathcal O}_{{\mathcal Y}_{\bullet}}$ with $sp_*{\mathcal O}_{]U_{\bullet}[_{{\mathcal Y}_{\bullet}}}$ we get the wanted isomorphism.\\

\section{The \v{C}ech double complex $B^{\bullet\bullet}$ and the Steenbrink double complex $A^{\bullet\bullet}$}
\label{doubledef}

\addtocounter{satz}{1}{{{}} \arabic{section}.\arabic{satz}}\newcounter{diff1}\newcounter{diff2}\setcounter{diff1}{\value{section}}\setcounter{diff2}{\value{satz}} For $e\ge 1$ and $\tau=\{\tau_1,\tau_2,\ldots,\tau_{e}\}\in S_{e}$ with $\tau_1<\tau_2<\ldots<\tau_{k+1}$ and $\gamma\in \tau$ define $pos(\gamma\in \tau)\in\mathbb{N}$ by $\tau_{pos(\gamma\in \tau)}=\gamma$. We continue to work with a fixed system $\{({\mathcal Z}_h,{\mathcal Y}_h)\}_{h\in H}$ as before. Let ${\bf C}\omega^q_{Y^i}=0$ if $q<0$ or $i\le 0$. For $s,t\in\mathbb{Z}$ we define the bidegree $(s,t)$-term of the {\it \v{C}ech double complex} $B^{\bullet\bullet}$ by $$B^{st}={\bf C}\omega^s_{Y^{t+1}}=\bigoplus_{\tau\in S_{t+1}}{\bf C}\omega^s_{Y^{t+1}_{\tau}}=\bigoplus_{\tau\in S_{t+1}}\frac{{\bf C}\widetilde{\omega}^s_{Y_{\tau}^{t+1}}}{{\bf C}\widetilde{\omega}^{s-1}_{Y_{\tau}^{t+1}}\wedge\theta}.$$The vertical differentials $$B^{q(i-1)}={{\bf C}}\omega^q_{Y^i}\stackrel{d}{\longrightarrow} B^{(q+1)(i-1)}={\bf C}\omega^{q+1}_{Y^{i}}$$are those of ${\bf C}\omega^{\bullet}_{Y^{i}}$; the horizontal differentials are $$B^{q(i-1)}={\bf C}\omega^q_{Y^i}=\bigoplus_{\sigma\in S_i}{\bf C}\omega_{Y_{\sigma}^i}^{q}\stackrel{\epsilon}{\longrightarrow}B^{qi}={\bf C}\omega_{Y^{i+1}}^{q}=\bigoplus_{\tau\in S_{i+1}}{\bf C}\omega_{Y_{\tau}^{i+1}}^{q}$$$$(\eta_{\sigma})_{\sigma\in S_i}\mapsto(\sum_{j\in\tau}(-1)^{pos(j\in\tau)+q+1}\eta_{(\tau-j)})_{\tau\in S_{i+1}}.$$Here $(\tau-j)$ means the element $\tau-\{j\}$ of $S_i$, and by $\eta_{(\tau-j)}$ we actually mean the image of $\eta_{(\tau-j)}$ under the obvious restriction map ${\bf C}\omega_{Y_{\tau-j}^i}^{q}\to {\bf C}\omega_{Y_{\tau}^{i+1}}^{q}$. We define the augmentation $${\bf C}\omega^{\bullet}_{Y}\stackrel{\epsilon_0}{\longrightarrow}{\bf C}\omega^{\bullet}_{Y^1}=B^{\bullet0}$$as the sum of the canonical restriction maps (similar to the above differentials $\epsilon$, but without alternating signs). It induces a morphism of complexes ${\bf C}\omega^{\bullet}_{Y}\to B^{\bullet}$ with $B^{\bullet}$ the total complex of $B^{\bullet\bullet}$.\\

\addtocounter{satz}{1}{{{}} \arabic{section}.\arabic{satz}} Copying \cite{mokr} we define the {\it Steenbrink double complex} $(A^{ij})_{i,j}$ by$$A^{ij}=\frac{{\bf C}\widetilde{\omega}^{i+j+1}_{Y}}{P_j{\bf C}\widetilde{\omega}^{i+j+1}_{Y}}$$ if $i\ge 0, j\ge 0$, and $=0$ else. The vertical differentials $A^{ij}\to A^{i+1,j}$ are induced by $(-1)^jd:{\bf C}\widetilde{\omega}^{i+j+1}_{{Y}}\to {\bf C}\widetilde{\omega}^{i+j+2}_{{Y}}$; the horizontal differentials $A^{ij}\to A^{i,j+1}$ are induced by the assignment $\omega\mapsto\omega\wedge\theta$. The augmentation $${\bf C}\omega^{\bullet}_{Y}\longrightarrow A^{\bullet0},\quad \quad \omega\mapsto\omega\wedge\theta$$ defines a morphism of complexes ${\bf C}\omega^{\bullet}_{Y}\to A^{\bullet}$ with $A^{\bullet}$ the total complex of $A^{\bullet\bullet}$.

\begin{pro}\label{vglabb} (1) ${\bf C}\omega^{\bullet}_{Y}\to B^{\bullet}$ is a quasiisomorphism.\\(2) ${\bf C}\omega^{\bullet}_{Y}\to A^{\bullet}$ is a quasiisomorphism.\\(3) There is a morphism of double complexes $\psi:A^{\bullet\bullet}\to B^{\bullet\bullet}$ compatible with the respective augmentations by ${\bf C}\omega^{\bullet}_{Y}$. It induces a quasiisomorphism$$\psi:A^{\bullet}\longrightarrow B^{\bullet}.$$In particular we may identify $H^{*}_{crys}(Y/T)_{\mathbb Q}=H^{*}(Y,{\bf C}\omega^{\bullet}_{{Y}})=H^{*}(Y,B^{\bullet})=H^{*}(Y,A^{\bullet})$. 
\end{pro}

{\sc Proof:} For sheaf complexes $K^{\bullet}$ on $U_{\bullet}=\{U_G\}_{\emptyset\ne G\subset H}$ and $\emptyset\ne G\subset H$ let $(K^{\bullet})_G$ be the component of $K^{\bullet}$ on $U_G$. For each $G$ the complex $({\bf C}\omega^{\bullet}_{Y})_{G}$, resp. $(B^{\bullet})_G$, resp. $(A^{\bullet})_G$ is quasiisomorphic with the corresponding complex formed with respect to a single admissible lift of $U_G$: this follows from Lemma \ref{dpring} and the Poincar\'{e} lemma for open polydisks (applied in the case of $(B^{\bullet})_G$, resp. $(A^{\bullet})_G$ to the subquotient complexes $(B^{\bullet,j})_G$, resp. $(A^{\bullet,j})_G$ for all $j$). Therefore we may assume for (1) and (2) that we are working with one global admissible lift $({\mathcal Z},{\mathcal Y})$ and drop $G$ from our notations. To show (1) it is enough to show that\begin{gather}0\longrightarrow {\bf C}\omega^q_Y\longrightarrow {\bf C}\omega^q_{Y^1}\longrightarrow {\bf C}\omega^q_{Y^{2}}\longrightarrow\ldots\tag{$*$}\end{gather}is exact for each $q$. Since ${\bf C}\omega^q_{Y^i}={\bf C}\omega^q_{Y}\otimes_{{\mathcal O}_{{\mathcal Y}}}{\mathcal O}_{{\mathcal Y}^i}$ and ${\bf C}\omega^q_{Y}$ is locally free over ${\mathcal O}_{{\mathcal Y}}$, this follows from the exactness (Chinese remainder theorem) of$$0\longrightarrow{\mathcal O}_{{\mathcal Y}}\longrightarrow{\mathcal O}_{{\mathcal Y}^1}\longrightarrow{\mathcal O}_{{\mathcal Y}^{2}}\longrightarrow\ldots.$$Also statement (2) is reduced to the exactness of $(*)$, literally as \cite{mokr} 3.16 is reduced to \cite{mokr} 3.15.1. Indeed, the proof in \cite{mokr} 3.16, although written for logarithmic {\it de Rham Witt} complexes, is in fact a completeley general argument valid for logarithmic de Rham complexes for any (relative) normal crossings divisor on a (relative) smooth (formal) scheme (note that since we are working with one global admissible lift we are not taking rigid analytic tubes of proper subschemes, just as in \cite{mokr} 3.16 one does not need to deal with divided power envelopes of proper subschemes). 

We turn to (3). The assertion on quasiisomorphy follows from (1) and (2) once $\psi$ with the stated ${\bf C}\omega^{\bullet}_Y$-compatibility is defined. To do this first note that for $\sigma\in S_k$ the residue map $Res_{\sigma}:\widetilde{\omega}_{{\mathcal K}_{\bullet}}^{q+k}\to\widetilde{\omega}_{{\mathcal K}_{\bullet}}^{q}\otimes{\mathcal O}_{{\mathcal Y}^k_{\bullet,\sigma}}$ described in \ref{logresi} induces a residue map $Res_{\rho}:{\bf C}\widetilde{\omega}_Y^{q+k}\to {\bf C}\widetilde{\omega}^q_{Y_{\sigma}^k}$, and by construction the latter vanishes on $P_{k-1}{\bf C}\widetilde{\omega}^{q+k}_Y$. In particular me may define$$A^{q(k-1)}=\frac{{\bf C}\widetilde{\omega}_Y^{q+k}}{P_{k-1}{\bf C}\widetilde{\omega}^{q+k}_Y}\stackrel{\psi}{\longrightarrow}B^{q(k-1)}=\bigoplus_{\sigma\in S_k}\frac{{\bf C}\widetilde{\omega}^q_{Y_{\sigma}^k}}{{\bf C}\widetilde{\omega}^{q-1}_{Y_{\sigma}^k}\wedge\theta}$$$$\eta\mapsto(\alpha_{q,k}.Res_{\sigma}(\eta))_{\sigma\in S_k}$$for $k\ge 1$, where we set $\alpha_{q,k}=-1$ if $(q-1,k)\in 2\mathbb{Z}\times 2\mathbb{Z}$, and $\alpha_{q,k}=1$ otherwise. We first verify that for the horizontal differentials $\wedge\theta$ of $A^{\bullet\bullet}$ and $\epsilon$ of $B^{\bullet\bullet}$ we have $\psi\circ(\wedge\theta)=\epsilon\circ\psi$. This can be done locally, in particular we may forget about distant components. So we may assume that there are $\{t_i\}_{1\le i\le l}\in{\mathcal O}_{{\mathcal K}_G}$, units $\{t_i\}_{l+1\le i\le g}\in{\mathcal O}^{\times}_{{\mathcal K}_G}$ (with $G\subset H$ fixed) and an order preserving identification between $R$ and $\{1,\ldots,l\}$ such that ${\mathcal Y}_{G,i}={\mathbb V}(t_i)$ for all $i\in R=\{1,\ldots,l\}$ and such that $\dlog(t_1),\ldots,\dlog(t_g)$ is an ${\mathcal O}_{{\mathcal K}_G}$-basis of $\widetilde{\omega}^1_{{\mathcal K}_G}$. We may more specifically assume (after multiplying one of $t_1,\ldots,t_l$ by an appropriate unit) that $$t=\prod_{1\le i\le l}t_i\in{\mathcal O}_{{\mathcal K}_G},$$ the image of the distinguished element $t\in {\mathcal O}_V$. We identify $R=\{1,\ldots,l\}$ with the set $S_1$ of subsets of $R=\{1,\ldots,l\}$ with precisely one element. Then $\theta=\sum_{1\le i\le l}\dlog(t_i)=\sum_{i\in S_1}\dlog(t_i)$. For $e\le g$ we denote by $\widetilde{S}_e$ the set of subsets of $\{1,\ldots,g\}$ with precisely $e$ elements. We write $$\dlog (t_{\nu})=\dlog(t_{\nu_e})\wedge\ldots\wedge\dlog(t_{\nu_1})$$ for $\nu\in \widetilde{S}_e$ with elements $\nu_e>\ldots>\nu_1$, and for $\gamma\in \nu$ we define $pos(\gamma\in \nu)\in\mathbb{N}$ by $\nu_{pos(\gamma\in \nu)}=\gamma$. By definition, for $\eta\in A^{q(k-1)}$ we have $$\psi(\eta\wedge\theta)=(\alpha_{q,k+1}\sum_{i\in S_1}Res_{\tau}(\eta\wedge\dlog(t_i)))_{\tau \in S_{k+1}}$$$$\epsilon(\psi(\eta))=(\alpha_{q,k}\sum_{i\in\tau}(-1)^{pos(i\in\tau)+q+1}Res_{\tau-i}\eta)_{\tau\in S_{k+1}}$$in $$B^{qk}={\bf C}\omega^q_{Y^{k+1}}=\bigoplus_{\tau\in S_{k+1}}\frac{{\bf C}\widetilde{\omega}^q_{Y^{k+1}_{\tau}}}{{\bf C}\widetilde{\omega}^{q-1}_{Y^{k+1}_{\tau}}\wedge\theta}.$$ Now $\eta$ can be written as a sum of elements of the type $\beta\wedge\dlog(t_{\alpha})$ with $\alpha\in \widetilde{S}_{q+k}$ and $\beta\in {\bf C}\widetilde{\omega}_Y^0$. Therefore we may assume $\eta=\beta\wedge\dlog(t_{\alpha})$ for some $\alpha\in \widetilde{S}_{q+k}$ and $\beta\in {\bf C}\widetilde{\omega}_Y^0$. We check equality of the $\tau$-components of the above expressions for fixed $\tau\in S_{k+1}$. In case $|\tau-(\alpha\cap\tau)|\ge2$ both of them vanish. Now consider the case $\tau=(\alpha\cap\tau)\cup i_0$ for some $i_0\notin\alpha$. We have $\eta=\beta'\wedge\dlog(t_{\alpha-(\alpha\cap\tau)})\wedge\dlog(t_{\alpha\cap\tau})$ with $\beta'=\beta$ or $\beta'=-\beta$. Then\begin{align}&\quad\quad\alpha_{q,k+1}\sum_{i\in S_1}Res_{\tau}(\eta\wedge\dlog(t_i))=\alpha_{q,k+1}Res_{\tau}(\eta\wedge\dlog(t_{i_0}))\notag\\&=\alpha_{q,k+1}(-1)^{pos(i_0\in\tau)+1}Res_{\tau}(\beta'\wedge\dlog(t_{\alpha-(\alpha\cap\tau)})\wedge\dlog(t_{\tau}))\notag\\&=\alpha_{q,k+1}(-1)^{pos(i_0\in\tau)+1}\beta'\wedge\dlog(t_{\alpha-(\alpha\cap\tau)})\notag\end{align}and$$\alpha_{q,k}\sum_{i\in\tau}(-1)^{pos(i\in\tau)+q+1}Res_{\tau-i}\eta=\alpha_{q,k}(-1)^{pos(i_0\in\tau)+q+1}\beta'\wedge\dlog(t_{\alpha-(\alpha\cap\tau)})$$and we see equality. Finally consider the case $\tau\subset\alpha$. We have $\eta=\beta'\wedge\dlog(t_{\alpha-\tau})\wedge\dlog(t_{\tau})$ with $\beta'=\beta$ or $\beta'=-\beta$. Then $$\alpha_{q,k+1}\sum_{i\in S_1}Res_{\tau}(\eta\wedge\dlog(t_i))=\alpha_{q,k+1}\sum_{i\in S_1-\alpha}Res_{\tau}(\eta\wedge\dlog(t_i))$$\begin{align}&=\alpha_{q,k+1}\sum_{i\in S_1-\alpha}(-1)^{k+1+pos(i\in i\cup(\alpha-\tau))+1}Res_{\tau}(\beta'\wedge\dlog(t_{i\cup(\alpha-\tau)})\wedge\dlog(t_{\tau}))\notag\\&=\alpha_{q,k+1}\sum_{i\in S_1-\alpha}(-1)^{k+pos(i\in i\cup(\alpha-\tau))}\beta'\wedge\dlog(t_{i\cup(\alpha-\tau)})\notag\end{align}and$$\alpha_{q,k}\sum_{i\in\tau}(-1)^{pos(i\in\tau)+q+1}Res_{\tau-i}\eta$$\begin{align}&=\alpha_{q,k}\sum_{i\in\tau}(-1)^{k+1+q+1}Res_{\tau-i}(\beta'\wedge\dlog(t_{\alpha-\tau})\wedge\dlog(t_i)\wedge\dlog(t_{\tau-i}))\notag\\&=\alpha_{q,k}\sum_{i\in\tau}(-1)^{k+q+pos(i\in i\cup(\alpha-\tau))+1}\beta'\wedge\dlog(t_{i\cup(\alpha-\tau)}).\notag\end{align}But in ${\bf C}\widetilde{\omega}^{q+1}_{Y_{\tau}^{k+1}}/({\bf C}\widetilde{\omega}^{q}_{Y_{\tau}^{k+1}}\wedge\theta)$ the element \begin{align}\beta'\wedge\dlog(t_{\alpha-\tau})\wedge\theta&=\sum_{i\in S_1}\beta'\wedge\dlog(t_{\alpha-\tau})\wedge\dlog(t_i)\notag\\{}&=\sum_{i\in S_1-(\alpha-\tau)}(-1)^{pos(i\in i\cup(\alpha-\tau))}\beta'\wedge\dlog(t_{i\cup(\alpha-\tau)})\end{align}vanishes, thus again we obtain the desired equality and the proof of $\psi\circ(\wedge\theta)=\epsilon\circ\psi$ is finished. The compatibility of $\psi$ with the vertical differentials $d$ is easy to check.\\Now we check  $\epsilon_0=\psi\circ(\wedge\theta)$ for the augmentation maps $\epsilon_0:{\bf C}\omega^q_Y\to B^{q0}$ and $(\wedge\theta):{\bf C}\omega^q_Y\to A^{q0}$. Let $\eta\in {\bf C}\omega^q_Y$. We work with a system of coordinates $\{t_i\}$ as above, and again we may assume $\eta=\beta\wedge\dlog(t_{\alpha})$ for some $\alpha\in \widetilde{S}_q$ and $\beta\in {\bf C}\omega_Y^0$. We find $$\psi(\eta\wedge\theta)=(Res_j\sum_{i\in S_1}\eta\wedge\dlog(t_i))_{j\in S_1}=(Res_j\sum_{i\in S_1-\alpha}\eta\wedge\dlog(t_i))_{j\in S_1}$$and we see that for $j\in S_1-\alpha$ the $j$-component is $\eta$, i.e. the $j$-component of $\epsilon_0(\eta)$. Now consider the respective $j$-components for $j\in\alpha$. We have $\eta=\beta'\wedge\dlog(t_{\sigma})\wedge\dlog(t_j)$ with  $\sigma=\alpha-j$ and $\beta'=\beta$ or $\beta'=-\beta$. Then the $j$-component of $\psi(\eta\wedge\theta)$ is$$Res_j\sum_{i\in S_1-\alpha}(-1)^{pos(i\in(\sigma\cup i))}\beta'\wedge\dlog(t_{\sigma\cup i})\wedge\dlog(t_j)$$$$=\sum_{i\in S_1-\alpha}(-1)^{pos(i\in(\sigma\cup i))}\beta'\wedge\dlog(t_{\sigma\cup i}).$$Since $$\beta'\wedge\dlog(t_{\sigma})\wedge\theta=\sum_{i\in S_1-\alpha}(-1)^{pos(i\in(\sigma\cup i))+1}\beta'\wedge\dlog(t_{\sigma\cup i})$$vanishes in ${\bf C}\omega_{Y^1}^q$, we find that also in this case the $j$-component is $\eta$, i.e. the $j$-component of $\epsilon_0(\eta)$.\\

\addtocounter{satz}{1}{{{}} \arabic{section}.\arabic{satz}} Both $B^{\bullet}$ and $A^{\bullet}$ are designed to compute $H^{*}_{crys}(Y/T)\otimes{\mathbb{Q}}$ in terms of cohomology groups of the components of $Y$. A first advantage of $B^{\bullet}$ against $A^{\bullet}$ is the following. There does not seem to exist an obvious pairing on $A^{\bullet}$, compatible with the differential in the usual sense, extending the cup product on ${\bf C}\omega_Y^{\bullet}$ which induces the Poincar\'{e}-duality in cohomology from \cite{hy}. By contrast, such a pairing does exist on $B^{\bullet}$. Namely, define$${\bf C}\omega^i_{Y^{j+1}}\otimes {\bf C}\omega^k_{Y^{s+1}}\longrightarrow {\bf C}\omega^{i+k}_{Y^{j+s+1}}$$by composing the restriction maps$${\bf C}\omega^i_{Y^{j+1}}\longrightarrow {\bf C}\omega^{i}_{Y^{j+s+1}}\quad\quad\mbox{and}\quad\quad {\bf C}\omega^k_{Y^{s+1}}\longrightarrow {\bf C}\omega^{k}_{Y^{j+s+1}}$$(for these we do {\it not} use alternating signs as we did in the definition of the differential $\epsilon$ of $B^{\bullet}$) with the cup product$${\bf C}\omega^i_{Y^{j+s+1}}\otimes {\bf C}\omega^k_{Y^{j+s+1}}\longrightarrow {\bf C}\omega^{i+k}_{Y^{j+s+1}}.$$

\addtocounter{satz}{1}{{{}} \arabic{section}.\arabic{satz}}\newcounter{cecmodn1}\newcounter{cecmodn2}\setcounter{cecmodn1}{\value{section}}\setcounter{cecmodn2}{\value{satz}} On the double complex $B^{\bullet\bullet}$ define the \v{C}ech filtration $F_C^{\bullet}B^{\bullet\bullet}$ by setting $F_C^rB^{ik}=B^{ik}$ if $r<k+1$, and $F_C^rB^{ik}=0$ if $r\ge k+1$. It gives rise to a spectral sequence\begin{gather}E_1^{pq}=H^q(Y,B^{\bullet p})\Longrightarrow H^{p+q}(Y,B^{\bullet})=H^{p+q}_{crys}(Y/T)_{\mathbb Q}.\tag*{$(C)_{T}$}\end{gather}We denote by $F_C^{\bullet}H^{*}_{crys}(Y/T)_{\mathbb Q}$ the induced filtration on $H^{*}_{crys}(Y/T)_{\mathbb Q}$.\\

\addtocounter{satz}{1}{{{}} \arabic{section}.\arabic{satz}} Let $\nu$ be the bihomogeneous endomorphism of bidegree (-1,1) of $A^{\bullet\bullet}$ such that $(-1)^{j+1}\nu$ is the natural projection $A^{i,j}\to A^{i-1,j+1}$. By \cite{mokr} 3.18 it induces the usual monodromy operator $N$ on 
$H_{crys}^*(Y/T)_{\mathbb Q}=H^*(Y,A^{\bullet})$. (In fact our definition of $\nu$ differs from the one written in \cite{mokr} 3.13 by the sign $(-1)^i$. It seems to us that to ensure that \cite{mokr} 3.18 holds, our definition is the correct one). $\nu$ {\it anti}commutes with the differential of $A^{\bullet\bullet}$. Note that the filtration $\nu^{\bullet}A^{\bullet\bullet}$ of $A^{\bullet\bullet}$ by the images of the iterated applications of $\nu$ is just the stupid horizontal filtration. Since the morphism $\psi$ from \ref{vglabb} sends $\nu^{r}A^{\bullet\bullet}$ to $F_C^{r}B^{\bullet\bullet}$ we see:

\begin{kor}\label{easyinkl} For all $r\ge0$,$$\bi(H^*(Y,\nu^rA^{\bullet})\longrightarrow H^{*}_{crys}(Y/T)_{\mathbb Q})\subset F_C^{r}H^{*}_{crys}(Y/T)_{\mathbb Q}$$ in $H^{*}_{crys}(Y/T)_{\mathbb Q}$. In particular $\bi{N^r}\subset F_C^{r}H^{*}_{crys}(Y/T)_{\mathbb Q}$ inside $H^{*}_{crys}(Y/T)_{\mathbb Q}$.
\end{kor}

\section{The \v{C}ech-Steenbrink tricomplex}
\label{trico}

\addtocounter{satz}{1}{{{}} \arabic{section}.\arabic{satz}} We now develop another tool to compare $B^{\bullet}$ with $A^{\bullet}$:  the {\it \v{C}ech-Steenbrink tricomplex} $C^{\bullet\bullet\bullet}$. For $k\ge 0, i\ge 0, j\ge 0$ its tridegree $(ijk)$-term is$$C^{ijk}=\frac{{\bf C}\widetilde{\omega}^{i+j+1}_{Y^{k+1}}}{P_j{\bf C}\widetilde{\omega}^{i+j+1}_{Y^{k+1}}}.$$For other triples $(i,j,k)$ we let $C^{ijk}=0$. The differentials $C^{ijk}\to C^{(i+1)jk}$ are those induced from $(-1)^jd:{\bf C}\widetilde{\omega}^{i+j+1}_{Y^{k+1}}\to {\bf C}\widetilde{\omega}^{i+j+2}_{Y^{k+1}}$. The differentials $C^{ijk}\to C^{i(j+1)k}$ are $\omega\mapsto\omega\wedge\theta$. The differentials $C^{ij(k-1)}\to C^{ijk}$ are$$C^{ij(k-1)}=\bigoplus_{\sigma\in S_k}\frac{{\bf C}\widetilde{\omega}^{i+j+1}_{Y^k_{\sigma}}}{P_j{\bf C}\widetilde{\omega}^{i+j+1}_{Y^k_{\sigma}}}\stackrel{\epsilon}{\longrightarrow}C^{ijk}=\bigoplus_{\tau\in S_{k+1}}\frac{{\bf C}\widetilde{\omega}^{i+j+1}_{Y^{k+1}_{\tau}}}{P_j{\bf C}\widetilde{\omega}^{i+j+1}_{Y^{k+1}_{\tau}}}$$$$(\eta_{\sigma})_{\sigma\in S_k}\mapsto(\sum_{\gamma\in\tau}(-1)^{pos(\gamma\in\tau)+i+1}\eta_{(\tau-\gamma)})_{\tau\in S_{k+1}}$$(here $(\tau-\gamma)$ means the element $\tau-\{\gamma\}$ of $S_k$, and by $\eta_{(\tau-\gamma)}$ we actually mean the image of $\eta_{(\tau-\gamma)}$ under the restriction map).\\
Fix $k\ge 1$. The differentials of the tricomplex $C^{\bullet\bullet\bullet}$ induce on its subsheaf $C^{\bullet\bullet(k-1)}$ a structure of double complex; we denote its associated total complex by $C^{\bullet(k-1)}$, {\it with the convention that $C^{ij(k-1)}$ is placed in degree $i+j$}. Similarly we can form the complex $B^{\bullet (k-1)}$ {\it with the convention that $B^{i(k-1)}$, when regarded as a position in $B^{\bullet (k-1)}$, is placed in degree $i$}. The augmentation $B^{\bullet(k-1)}\to C^{\bullet 0(k-1)}$ is defined by $\omega\mapsto\omega\wedge\theta$.

\begin{pro}\label{diqui} This induces a quasiisomorphism $B^{\bullet(k-1)}\to C^{\bullet(k-1)}$.
\end{pro}

{\sc Proof:} It is enough to show that for all $q\ge 0$ the sequence$$0\longrightarrow {\bf C}\omega^q_{Y^k}\longrightarrow C^{q,0,(k-1)}\longrightarrow C^{q,1,(k-1)}\longrightarrow C^{q,2,(k-1)}\longrightarrow\ldots$$is exact, in other words that$${\bf C}\widetilde{\omega}^{q-1}_{Y^k}\stackrel{\wedge\theta}{\longrightarrow}{\bf C}\widetilde{\omega}^q_{Y^k}\stackrel{\wedge\theta}{\longrightarrow}\frac{{\bf C}\widetilde{\omega}^{q+1}_{Y^k}}{P_0{\bf C}\widetilde{\omega}^{q+1}_{Y^k}}\stackrel{\wedge\theta}{\longrightarrow}\frac{{\bf C}\widetilde{\omega}^{q+2}_{Y^k}}{P_1{\bf C}\widetilde{\omega}^{q+2}_{Y^k}}\stackrel{\wedge\theta}{\longrightarrow}\ldots$$is exact. As in \cite{mokr} 3.15 this is reduced to proving that
$$0\longrightarrow\Gr_0{\bf C}\widetilde{\omega}^{\bullet}_{Y^k}\stackrel{\wedge\theta}{\longrightarrow}\Gr_1{\bf C}\widetilde{\omega}^{\bullet}_{Y^k}[1]\stackrel{\wedge\theta}{\longrightarrow}\Gr_2{\bf C}\widetilde{\omega}^{\bullet}_{Y^k}[2]\stackrel{\wedge\theta}{\longrightarrow}\ldots$$is an exact sequence of complexes. In view of \ref{corand} it is enough to show that for each component $Y^k_s$ of $Y^k$ the complex$$0\longrightarrow {\bf C}\Omega^{\bullet}_{Y_s^k}\longrightarrow \bigoplus_{t\in S_1}{\bf C}\Omega^{\bullet}_{Y_t^1\cap Y_s^k}\longrightarrow\bigoplus_{t\in S_2}{\bf C}\Omega^{\bullet}_{Y_t^2\cap Y_s^k}\longrightarrow\ldots$$ is exact. The differentials in this complex are the sums of all restriction maps (with alternating signs) between the respective summands which make sense. It is then enough to show that for arbitrary $q$ the complex$$M^{\bullet}=[{\bf C}\Omega^q_{Y_s^k}\longrightarrow \bigoplus_{t\in S_1}{\bf C}\Omega^q_{Y_t^1\cap Y_s^k}\longrightarrow\bigoplus_{t\in S_2}{\bf C}\Omega^q_{Y_t^2\cap Y_s^k}\longrightarrow\ldots]$$is exact. This is the total complex of the following double complex $M^{\bullet\bullet}$. Let $J=\{u\in R|\, Y_s^k\subset Y_u\}$ and $I=R-J$. Let $$M^{ij}=\bigoplus_{{t\in S_{i+j}}\atop{|t\cap I|=j}}{\bf C}\Omega^q_{Y_t^{i+j}\cap Y_s^k}$$for $i\ge 0$ and $j\ge 0$, and $M^{ij}=0$ otherwise. The differentials $M^{ij}\to M^{i+1,j}$ and $M^{ij}\to M^{i,j+1}$ are the sums of all restriction maps (with alternating signs) between the respective summands which make sense. Now it is enough to show that for each fixed $j$ the complex $M^{\bullet j}$ is exact. We have a direct sum decomposition of complexes$$M^{\bullet j}=\bigoplus_{{r\in S_j}\atop{r\subset I}}M_r^{\bullet}$$
$$M_r^i=\bigoplus_{{t\in S_i}\atop{t\subset J}}{\bf C}\Omega^q_{Y_r^{j}\cap Y_s^k}$$where the differentials $M_r^i\to M_r^{i+1}$ are the sums of identity maps (with appropriate alternating signs) between the respective summands. Since the complexes $M_r^{\bullet}$ are exact, so is $M^{\bullet j}$.\\We remark that as in \ref{vglabb} we could also have reduced to the case where we are dealing with a single global admissible lift. At least in that case the complexes $M^{i \bullet}$ are exact in positive degrees (with $M^{i0}$ in degree $0$). We do not need this.\\  

\addtocounter{satz}{1}{{{}} \arabic{section}.\arabic{satz}} Let $C^{\bullet}$ be the total complex of the tricomplex $C^{\bullet\bullet\bullet}$. The augmentation $$A^{\bullet\bullet}\longrightarrow C^{\bullet\bullet 0}$$ is defined as the sum of the canonical restriction maps (similar to the differentials $C^{ij(k-1)}\to C^{ijk}$, but without alternating signs). By \ref{vglabb} and \ref{diqui} we have:

\begin{kor}\label{allquis} The augmentations $A^{\bullet\bullet}\to C^{\bullet\bullet 0}$ and $B^{\bullet\bullet}\to C^{\bullet 0\bullet}$ induce quasiisomorphisms $A^{\bullet}\to C^{\bullet}$ and $B^{\bullet}\to C^{\bullet}$. In particular we may identify $H^*_{crys}(Y/T)_{\mathbb Q}=H^*(Y,{\bf C}\omega^{\bullet}_Y)=H^*(Y,A^{\bullet})=H^*(Y,B^{\bullet})=H^*(Y,C^{\bullet})$.
\end{kor}

\addtocounter{satz}{1}{{{}} \arabic{section}.\arabic{satz}} Let $\nu$ be the trihomogeneous endomorphism of tridegree (-1,1,0) of $C^{\bullet\bullet\bullet}$ such that $(-1)^{j+1}\nu$ is the natural projection $C^{ij(k-1)}\to C^{(i-1)(j+1)(k-1)}$. It is clear that the resolution $A^{\bullet}\to C^{\bullet}$ is compatible with the endomorphisms $\nu$ on source and target, thus $\nu$ on $C^{\bullet\bullet\bullet}$ induces the same endomorphism $N$ in cohomology. For fixed $k\ge1$ we also denote by $N$ the operator on $H^*(Y,{\bf C}\omega_{Y^k}^{\bullet})$ induced by the endomorphism $\nu:C^{\bullet\bullet(k-1)}\to C^{\bullet\bullet(k-1)}$ via the canonical isomorphism $H^*(Y,C^{\bullet(k-1)})=H^*(Y,B^{\bullet(k-1)})\cong H^*(Y,{\bf C}\omega_{Y^k}^{\bullet})$.

\begin{satz}\label{mokomptri} For any $k\ge 1$ we have $N=0$ on $H^*(Y,{\bf C}\omega_{Y^k}^{\bullet})$. 
\end{satz}

{\sc Proof:} $C^{\bullet\bullet(k-1)}$ is the direct sum, over all $\sigma\in S_k$, of the double complexes$$C^{\bullet\bullet(k-1)}_{\sigma}=\bigoplus_{i,j}C^{ij(k-1)}_{\sigma}\quad\mbox{ with }\quad C^{ij(k-1)}_{\sigma}=(\frac{{\bf C}\widetilde{\omega}^{i+j+1}_{Y_{\sigma}^k}}{P_j{\bf C}\widetilde{\omega}^{i+j+1}_{Y_{\sigma}^k}})_{i,j}$$and $\nu$ is given by endomorphisms $\nu$ on each of these summands. We will show that on each summand $\nu$ is homotopic to zero. Fix $\sigma\in S_k$, choose an auxiliary $\gamma\in\sigma$ and define$$C^{ij(k-1)}_{\sigma}\stackrel{h^{ij}}{\longrightarrow}C^{(i-1)j(k-1)}_{\sigma}$$$$\omega\mapsto (-1)^{j+1}Res_{\gamma}\omega$$(this depends on $\gamma$). We claim \begin{gather}\nu(\omega)=h^{ij}(\omega)\wedge\theta-h^{i(j+1)}(\omega\wedge\theta)\tag{$1$}\end{gather} for $\omega\in C^{ij(k-1)}_{\sigma}$. This can be verified locally, in particular we may forget about distant components. So we may assume that there are $\{t_i\}_{1\le i\le l}\in{\mathcal O}_{{\mathcal K}_G}$, units $\{t_i\}_{l+1\le i\le g}\in{\mathcal O}^{\times}_{{\mathcal K}_G}$ (with $G\subset H$ fixed) and an order preserving identification between $R$ and $\{1,\ldots,l\}$ such that ${\mathcal Y}_{G,i}={\mathbb V}(t_i)$ for all $i\in R=\{1,\ldots,l\}$ and such that $\dlog(t_1),\ldots,\dlog(t_g)$ is an ${\mathcal O}_{{\mathcal K}_G}$-basis of $\widetilde{\omega}^1_{{\mathcal K}_G}$. We may more specifically assume (after multiplying one of $t_1,\ldots,t_l$ by an appropriate unit) that $$t=\prod_{1\le i\le l}t_i\in{\mathcal O}_{{\mathcal K}_G},$$ the image of the distinguished element $t\in {\mathcal O}_V$. We identify $R=\{1,\ldots,l\}$ with the set $S_1$ of subsets of $R=\{1,\ldots,l\}$ with precisely one element. Then $\theta=\sum_{1\le i\le l}\dlog(t_i)=\sum_{i\in S_1}\dlog(t_i)$. For $e\le g$ we denote by $\widetilde{S}_e$ the set of subsets of $\{1,\ldots,g\}$ with precisely $e$ elements. We write $$\dlog (t_{\nu})=\dlog(t_{\nu_e})\wedge\ldots\wedge\dlog(t_{\nu_1})$$ for $\nu\in \widetilde{S}_e$ with elements $\nu_e>\ldots>\nu_1$, and for $\gamma\in \nu$ we define $pos(\gamma\in \nu)\in\mathbb{N}$ by $\nu_{pos(\gamma\in \nu)}=\gamma$. We may assume that $\omega$ is represented by $\beta\dlog(t_{\rho})$ for some $\beta\in {\bf C}\omega_{Y_{\sigma}^k}^0$ and $\rho\in \widetilde{S}_{i+j+1}$. In case $\gamma\in\rho$ we find $h^{ij}(\omega)=(-1)^{pos(\gamma\in\rho)+j}\beta\dlog(t_{\rho-\gamma})$, hence \begin{gather}h^{ij}(\omega)\wedge\theta=\sum_{\alpha\in S_1-(\rho-\gamma)}(-1)^{pos(\gamma\in\rho)+j}\beta\dlog(t_{\rho-\gamma})\wedge\dlog(t_{\alpha}).\tag{$i$}\end{gather}On the other hand $\omega\wedge\theta=\sum_{\alpha\in S_1-\rho}\beta\dlog(t_{\rho})\wedge\dlog(t_{\alpha})$ and hence \begin{gather}h^{i(j+1)}(\omega\wedge\theta)=\sum_{\alpha\in S_1-\rho}(-1)^{pos(\gamma\in\rho)+j}\beta\dlog(t_{\rho-\gamma})\wedge\dlog(t_{\alpha}).\tag{$ii$}\end{gather} The difference $(i)-(ii)$ is $(-1)^{j+1}\beta\dlog(t_{\rho})=\nu(\omega)$. In case $\gamma\notin\rho$ we see $h^{ij}(\omega)=0$, while as before $\omega\wedge\theta=\sum_{\alpha\in S_1-\rho}\beta\dlog(t_{\rho})\wedge\dlog(t_{\nu})$ and thus $h^{i(j+1)}(\omega\wedge\theta)=(-1)^{j+1}\beta\dlog(t_{\rho})=\nu(\omega)$, so $(1)$ is proved. One also has \begin{gather} d(h^{ij}(\omega))-h^{(i+1)j}(d\omega)=0\tag{$2$}\end{gather} where $d$ denotes the differential $C^{\bullet\bullet(k-1)}_{\sigma}\to C^{(\bullet+1)\bullet(k-1)}_{\sigma}$. Now $(1)$ and $(2)$ together imply that on the complex $C^{\bullet(k-1)}_{\sigma}$, if $\widetilde{d}$ denotes its differential, we have $\nu=\widetilde{d}\circ h - h\circ\widetilde{d}$. This tells us that the endomorphisms $\nu$ of the graded algebra $C^{\bullet(k-1)}$, which anticommutes with $\widetilde{d}$, induces the zero map in cohomology.\\

\addtocounter{satz}{1}{{{}} \arabic{section}.\arabic{satz}} On $C^{\bullet\bullet\bullet}$ define the \v{C}ech filtration $F_C^{\bullet}C^{\bullet\bullet\bullet}$ by setting $F_C^rC^{ijk}=C^{ijk}$ if $r<k+1$, and $F_C^rC^{ijk}=0$ if $r\ge k+1$. Then \ref{diqui} says that $B^{\bullet}\to C^{\bullet}$ is a filtered quasiisomorphism with respect to the respective \v{C}ech filtrations. {\it In general}, for filtrations (on various complexes) denoted $P_{\bullet}$, resp. $F_C^{\bullet}$, we will use the notation $\Gr_{\bullet}$, resp. $\Gr^{\bullet}_C$, for the associated graded object.

\begin{kor}\label{inkl} $N(F_C^kH^*_{crys}(Y/T)_{\mathbb{Q}})\subset F_C^{k+1}H^*_{crys}(Y/T)_{\mathbb{Q}}$ for $k\ge0$ (sharpening \ref{easyinkl}).
\end{kor}

{\sc Proof:} It is enough to show that for $k\ge 0$, the image of $$H^*(Y,F^k_C B^{\bullet})\stackrel{H^*(\nu)}{\longrightarrow}H^*(Y,F^k_CB^{\bullet})$$is contained in the image of the natural map $H^*(Y,F^{k+1}_CB^{\bullet})\stackrel{\iota}{\to} H^*(Y,F^k_CB^{\bullet})$. In view of the exact sequence$$H^*(Y,F^{k+1}_CB^{\bullet})\stackrel{\iota}{\longrightarrow} H^*(Y,F^k_CB^{\bullet})\stackrel{pr}{\longrightarrow}H^*(Y,\Gr^k_CB^{\bullet})$$ this is equivalent with showing that the image of $H^*(\nu)$ is contained in $\ke(pr)$. This follows from \ref{mokomptri} using the canonical quasiisomorphism $B^{\bullet k}[k]\cong \Gr^k_CB^{\bullet}$.

\section{Monodromy and weight filtration via the \v{C}ech complex $B^{\bullet}$}
\label{rigdesc}

\addtocounter{satz}{1}{{{}} \arabic{section}.\arabic{satz}}\newcounter{wess1}\newcounter{wess2}\setcounter{wess1}{\value{section}}\setcounter{wess2}{\value{satz}} The endomorphisms $\Phi:=F$ on the complexes ${\bf C}{\omega}_{Y^i}^{\bullet}$ induce a Frobenius endomorphism $\Phi$ on $B^{\bullet}$. The Frobenius endomorphism $\Phi$ on $A^{\bullet}$ is defined by $\Phi:=p^{i}F$ on $A^{ij}$. On $A^{\bullet}$ the weight filtration $P_{\bullet}A^{\bullet}$ is defined by$$P_kA^{\bullet}=\bigoplus_{i\ge0, j\ge0}P_kA^{ij}\quad\quad\quad P_kA^{ij}=\frac{P_{2j+k+1}{\bf C}\widetilde{\omega}^{i+j+1}_Y}{P_j{\bf C}\widetilde{\omega}^{i+j+1}_Y}.$$Denote by $\Gr_{\bullet}A^{\bullet}$ the associated graded with respect to $P_{\bullet}A^{\bullet}$. We have $$\Gr_kA^{\bullet}\cong\bigoplus_{{j\ge0}\atop{j\ge-k}}{\bf C}\Omega^{\bullet}_{Y^{2j+k+1}}[-2j-k](-j-k)$$via residue maps (see \cite{mokr} 3.22). Passing to the limit $n\to\infty$ and tensoring with $\mathbb{Q}$ we get the weight spectral sequence$$E_1^{-k,i+k}=\bigoplus_{{j\ge0}\atop{j\ge-k}}H^{i-2j-k}(Y,{\bf C}\Omega^{\bullet}_{Y^{2j+k+1}})(-j-k)\Longrightarrow H^i_{crys}(Y/T)_{\mathbb{Q}}.$$
It gives rise to the {\it weight filtration} on $H^*_{crys}(Y/T)_{\mathbb{Q}}$. On the other hand, the operator $N$ on $H^*_{crys}(Y/T)_{\mathbb{Q}}$ gives rise to the {\it monodromy filtration} on $H^*_{crys}(Y/T)_{\mathbb{Q}}$, defined in \cite{mokr} 3.26. Consider the statement:\begin{description}\item[(MW)] The weight spectral sequence degenerates in $E_2$, and the weight filtration coincides with the monodromy filtration on $H^*_{crys}(Y/T)_{\mathbb{Q}}.$\end{description}
The standard conjecture (\cite{mokr} 3.24, 3.27) in this context states that (MW) holds whenever $Y$ is the reduction of a projective semistable $A$-scheme as in section \ref{lifsec}. Nakkajima \cite{nakadeg} proved the degeneration in $E_2$ unconditionally.\\We see that $N$ and the weight filtration on $H^*_{crys}(Y/T)_{\mathbb{Q}}$ can be obtained from corresponding structures on $A^{\bullet}$. The following Theorems \ref{viawb} and \ref{fiwbmo} tell us that we can recover $N$ and the weight filtration also from structures on $B^{\bullet}$. 

\begin{satz}\label{viawb} For $i\ge 0$ the $i$-fold iterated monodromy operator $N^i$ on $H^*_{crys}(Y/T)_{\mathbb Q}=H^*(Y,{\bf C}\omega^{\bullet}_Y)=H^*(Y,B^{\bullet})$ is induced by a composite of complex morphisms$${\bf C}\omega^{\bullet}_{Y}\stackrel{\rho_i}{\longrightarrow} F_C^iB^{\bullet}\stackrel{\iota}{\longrightarrow}B^{\bullet},$$where $\rho_i$ is given by explicit residue maps and $\iota$ is the inclusion.
\end{satz}

{\sc Proof:} Consider the quasiisomorphisms$${\bf C}\omega^{\bullet}_{Y}\stackrel{\wedge\theta}{\longrightarrow}A^{\bullet}\quad\quad\mbox{and}\quad\quad A^{\bullet}\stackrel{\psi}{\longrightarrow}B^{\bullet}.$$The operator $N^i$ on $H^{*}_{crys}(Y/T)_{\mathbb{Q}}$ is the map in cohomology induced by the composite $\psi\circ\nu^i\circ(\wedge\theta)$. This is the following map of complexes $\rho_{i}:{\bf C}\omega^{\bullet}_{Y}\to F_C^iB^{\bullet}$: An element $\eta\in {\bf C}\omega_Y^q$ is sent to$$\epsilon(Res_{\sigma}(\eta\wedge\theta))_{\sigma\in S_{i+1}}\in \bigoplus_{\sigma\in S_{i+1}}{\bf C}\omega_{Y^{i+1}_{\sigma}}^{q-i}=B^{(q-i)i}$$with the sign $\epsilon=\alpha_{(q-i),(i+1)}\prod_{j=0}^{i-1}(-1)^{j+1}\in\{\pm1\}$ (note $\eta\wedge\theta\in {\bf C}\widetilde{\omega}^{q+1}_Y$).\\

\addtocounter{satz}{1}{{{}} \arabic{section}.\arabic{satz}} On $B^{\bullet\bullet}$ define the {\it canonical filtration} by setting$$P^{can}_lB^{i(k-1)}=\frac{P_{l-1}{\bf C}\widetilde{\omega}^i_{Y^k}}{({\bf C}\widetilde{\omega}^{i-1}_{Y^k}\wedge\theta)\cap P_{l-1}{\bf C}\widetilde{\omega}^i_{Y^k}}$$for $l\ge 0$. It induces a filtration $P^{can}_{\bullet}B^{\bullet}$ on $B^{\bullet}$. On $A^{\bullet\bullet}$ define the {\it canonical filtration} as the kernel filtration for $\nu$, i.e. $$P^{can}_lA^{\bullet\bullet}=\ke(\nu^l:A^{\bullet\bullet}\to A^{\bullet\bullet}),$$ or equivalently $P^{can}_lA^{ij}=P_{l-j-1}A^{ij}$ for $l\ge 0$. For any filtration denoted $P_{\bullet}^{can}$ (here and on various other complexes below) we write $\Gr_{\bullet}^{can}$ for the associated graded object.\\

\begin{lem}\label{canqui} The map $\psi:A^{\bullet}\to B^{\bullet}$ from \ref{vglabb} is a filtered quasiisomorphism with respect to canonical filtrations.
\end{lem}

{\sc Proof:} On ${\bf C}\omega^{\bullet}_Y$ define the canonical filtration by setting $$P^{can}_l{\bf C}\omega^i_{Y}=\frac{P_{l-1}{\bf C}\widetilde{\omega}^i_{Y}}{({\bf C}\widetilde{\omega}^{i-1}_{Y}\wedge\theta)\cap P_{l-1}{\bf C}\widetilde{\omega}^i_{Y}}=\frac{P_{l-1}{\bf C}\widetilde{\omega}^i_{Y}}{P_{l-2}{\bf C}\widetilde{\omega}^{i-1}_{Y}\wedge\theta}.$$for $l\ge 0$. The second equality in this definition can be justified by induction on $l$, noting that the natural surjection$$\frac{P_{l-1}{\bf C}\widetilde{\omega}^i_{Y}}{(P_{l-2}{\bf C}\widetilde{\omega}^{i-1}_{Y}\wedge\theta)+(P_{l-2}{\bf C}\widetilde{\omega}^{i}_{Y})}\longrightarrow\frac{P_{l-1}{\bf C}\widetilde{\omega}^i_{Y}}{(({\bf C}\widetilde{\omega}^{i-1}_{Y}\wedge\theta)\cap P_{l-1}{\bf C}\widetilde{\omega}^i_{Y})+(P_{l-2}{\bf C}\widetilde{\omega}^{i}_{Y})}$$is also injective because its composition with the map$$\frac{P_{l-1}{\bf C}\widetilde{\omega}^i_{Y}}{(({\bf C}\widetilde{\omega}^{i-1}_{Y}\wedge\theta)\cap P_{l-1}{\bf C}\widetilde{\omega}^i_{Y})+(P_{l-2}{\bf C}\widetilde{\omega}^{i}_{Y})}\stackrel{\wedge\theta}{\longrightarrow}\frac{P_{l}{\bf C}\widetilde{\omega}^{i+1}_{Y}}{P_{l-1}{\bf C}\widetilde{\omega}^{i+1}_{Y}}$$is injective by \cite{mokr} 3.15. Using the second defining expression for $P^{can}_l{\bf C}\omega^{\bullet}_{Y}$ we get from \cite{mokr} 3.15 that $\wedge\theta:{\bf C}\omega^{\bullet}_Y\to A^{\bullet}$ is a filtered resolution w.r.t. canonical filtrations. Using the first one we get that the same for $\epsilon_0:{\bf C}\omega^\bullet_Y\to B^{\bullet}$. The lemma follows.\\

{\it For the rest of this section we assume that $k$ is finite and $Y$ is proper and satisfies (MW).} We fix a cohomology degree $*$ and simply write $H^*$ for $H^*_{crys}(Y/T)_{\mathbb{Q}}$.

\begin{pro}\label{monweicon} (Chiarellotto, \cite{chia}) (1) Let $r\in \mathbb{N}$. The sequence $$H^*(Y,P^{can}_rA^{\bullet})\longrightarrow H^*\stackrel{N^r}{\longrightarrow}H^*\longrightarrow H^*(Y,\koke(\nu^r))$$is exact. In particular we have $$\bi(H^*(Y,\nu^rA^{\bullet})\longrightarrow H^*)=\bi(H^*\stackrel{N^r}{\longrightarrow}H^*).$$(2) We have$$H_{rig}^*(Y)=H_{rig}^*(Y/K_0)=H^*(Y,P^{can}_1A^{\bullet}).$$ 
\end{pro}

\addtocounter{satz}{1}{{{}} \arabic{section}.\arabic{satz}} Define the {\it \v{C}ech filtration} on $H_{rig}^*(Y)$ by setting$$F_C^sH_{rig}^*(Y)=\bi(H^*(Y,P^{can}_1F^s_CA^{\bullet})\longrightarrow H_{rig}^*(Y))$$for $s\ge0$, where we write $F^s_CA^{\bullet}=\nu^sA^{\bullet}$. It deserves indeed its name: It is the filtration arising from the spectral sequence 
\begin{gather}E_1^{pq}=H_{rig}^q(Y^{p+1})\Longrightarrow H_{rig}^{p+q}(Y)\tag*{$(C)_{rig}$}\end{gather}(we do not need this fact). Let $H^*_{rig}(Y)\stackrel{\iota}{\to}H^*$ be the canonical map. By Chiarellotto's result \cite{chia}, since we assume (MW), we have $\bi(\iota)=\ke(N)$. 

\begin{pro}\label{hrigana} For all $r\ge 0$ we have $$\iota(F_C^rH_{rig}^*(Y))=\ke(N)\cap\bi(N^r).$$
\end{pro}

{\sc Proof:} This follows formally from the fact that the filtrations $F^{\bullet}_CH^*_{rig}(Y)$ on $H^*_{rig}(Y)$ and $\ke(N)\cap\bi(N^{\bullet})$ on $\bi(\iota)=\ke(N)$ are in fact weight filtrations for the Frobenius actions, with the same strictly monoton increasing sequences of weights on the graded pieces. Indeed, we have $$P^{can}_1\Gr_C^rA^{\bullet}\cong\Gr_{-r}A^{\bullet r}[-r]\cong \Gr_{r+1}{\bf C}\widetilde{\omega}^{\bullet}_Y[1](r+1)\cong {\bf C}\Omega^{\bullet-r}_{Y^{r+1}}$$(as in \cite{mokr} 3.22(2)). Therefore $H^*(Y,P^{can}_1\Gr_C^rA^{\bullet})$ is pure of weight $*-r$ in view of \cite{chile}, hence also its subquotient $\Gr_C^rH_{rig}^*(Y)$. On the other hand, by (MW) the subquotient $(\ke(N)\cap\bi(N^r))/(\ke(N)\cap\bi(N^{r+1}))$ of $H^*$ is pure of weight $*-r$.\\ 

\begin{satz}\label{fiwbmo} Denote by $G_{\bullet}$ the convolution of the filtrations $P^{can}_{\bullet}$ and $F_C^{\bullet}$ on $B^{\bullet}$, i.e. the filtration defined by$$G_k=\sum_{i\ge 0}P_{i+1}^{can}F^{i-k}_CB^{\bullet}.$$Then $G_{\bullet}$ induces the monodromy filtration on $H^*_{crys}(Y/T)_{\mathbb{Q}}=H^*(Y,B^{\bullet})$.
\end{satz}

{\sc Proof:} We will show that $H^*(Y,\Gr^{can}_k\Gr^r_CB^{\bullet})$ is pure of weight $*-r+k-1$. This then implies that the filtration induced by $G_{\bullet}$ is the weight filtration for Frobenius on $H^*_{crys}(Y/T)_{\mathbb{Q}}$, i.e. the uniquely determined Frobenius stable filtration whose subquotients are pure and of mutually different weights; but this is the monodromy filtration, by (MW).\\Define the double complex$$\widetilde{B}^{\bullet\bullet}=({\bf C}\widetilde{\omega}^q_{Y^{i+1}})_{q,i}$$with differentials analoguous to those in \arabic{diff1}.\arabic{diff2}. Let $\widetilde{B}^{\bullet}$, $\widetilde{B}^{\bullet}$ be the associated total complexes. As in \arabic{cecmodn1}.\arabic{cecmodn2} define the \v{C}ech filtration $F_C^{\bullet}\widetilde{B}^{\bullet}$ by setting $F_C^r\widetilde{B}^{i(k-1)}=\widetilde{B}^{i(k-1)}$ if $r<k$, and $F_C^r\widetilde{B}^{i(k-1)}=0$ if $r\ge k$. Define the canonical filtration by setting $P^{can}_l\widetilde{B}^{i(k-1)}=P_{l-1}{\bf C}\widetilde{\omega}^i_{Y^k}$ for $l\ge 0$.\\
{\it Claim: $H^*(Y,\Gr^{can}_k\Gr^r_C\widetilde{B}^{\bullet})$ is pure of weight $*-r+k-1$, for all $k\ge 1$. For all $k\ge 2$ also $H^*(Y,\Gr^{can}_kF^r_C\widetilde{B}^{\bullet})$ is pure of weight $*-r+k-1$.}\\
The assertion in case $k=1$ follows from the proof of \ref{hrigana}, thanks to the isomorphism $\psi:P^{can}_1A^{\bullet}\cong P^{can}_1B^{\bullet}=P^{can}_1\widetilde{B}^{\bullet}$. Now let $k\ge 2$. We have$$\Gr_k^{can}\Gr_C^s\widetilde{B}^{\bullet}=\Gr_{k-1}{\bf C}\widetilde{\omega}^{\bullet}_{Y^{s+1}}[-s]\cong {\bf C}\Omega^{\bullet}_{Y^{s+1}\cap Y^{k-1}}[-s-k+1](-k+1).$$So $H^*(Y,\Gr^{can}_k\Gr^r_C\widetilde{B}^{\bullet})$ is pure of weight $*-s-k+1-2(-k+1)=*-r+k-1$. It follows that $H^*(Y,\Gr^{can}_kF^r_C\widetilde{B}^{\bullet})$ is mixed with weights at most $*-r+k-1$, and $H^*(Y,\Gr^{can}_k(\widetilde{B}^{\bullet}/F^{r}_C\widetilde{B}^{\bullet}))$ is mixed with weights at least $*-r+k$. On the other hand, $\Gr^{can}_k \widetilde{B}^{\bullet}$ is quasiisomorphic with $\Gr_{k+1}{\bf C}\widetilde{\omega}^{\bullet}_Y$ via the canonical augmentation ${\bf C}\widetilde{\omega}_Y^{\bullet}\to \widetilde{B}^{\bullet 0}$. But $\Gr_{k+1}{\bf C}\widetilde{\omega}^{\bullet}_Y\cong {\bf C}\Omega^{\bullet}_{Y^{j}}[-j](-j)$, thus $H^*(Y,\Gr_k^{can}{\bf C}\widetilde{\omega}^{\bullet}_Y)$ is pure of weight $*+k-1$. Now the long cohomology sequence associated to$$0\longrightarrow\Gr^{can}_kF^r_C\widetilde{B}^{\bullet}\longrightarrow\Gr^{can}_k \widetilde{B}^{\bullet}\longrightarrow\Gr^{can}_k(\widetilde{B}^{\bullet}/F^{r}_C\widetilde{B}^{\bullet})\longrightarrow0$$gives the second statement of the claim.\\
To prove that $H^*(Y,\Gr^{can}_k\Gr^r_CB^{\bullet})$ is pure of weight $*-r+k-1$ it is now enough to show that for all $k\ge 1$ the natural map$$H^*(Y,\Gr^{can}_k\Gr^r_C\widetilde{B}^{\bullet})\longrightarrow H^*(Y,\Gr^{can}_k\Gr^r_CB^{\bullet})$$is surjective. For this we need to show that the connecting maps in cohomology associated with the short exact sequence $$0\longrightarrow\Gr_{k-1}^{can}\Gr_C^rB^{\bullet-1}\stackrel{\wedge\theta}{\longrightarrow}\Gr_{k}^{can}\Gr_C^r\widetilde{B}^{\bullet}\longrightarrow\Gr_{k}^{can}\Gr_C^rB^{\bullet}\longrightarrow0$$are zero. This follows from \ref{mokomptri}, because as in \cite{mokr} 3.18 one sees that these connecting maps are induced by the endomorphism $\nu$ of $\Gr_C^rC^{\bullet}$.\\

\section{\v{C}ech filtration versus $\bi(N^{\bullet})$ filtration}\label{seccec}
\label{gleich} 

We give a sufficient criterion for the equality $\bi(N^r)=F_C^rH^*$, but we also give an example with {\it $\bi(N^r)\ne F_C^rH^*$ although (MW) holds true}. 

\begin{lem}\label{nillem} Let $V$ be an abelian group, $m, n\in\mathbb{N}$, $n\le m$, let $(0)=F^m\subset\ldots\subset F^1\subset F^0=V$ be a descending filtration, and let $N\in\ho(V,V)$ such that $N(F^{i-1})\subset F^i$ and $\ke(N)\cap N^iV=\ke(N)\cap F^i$ for all $i\ge n$. Then $N^iV=F^i$ for all $i\ge n$.
\end{lem}

{\sc Proof:} Descending induction on $i$: Let $x\in F^{i-1}$. Then $Nx=N^iy$ for some $y$ by induction hypothesis, thus $x-N^{i-1}y\in \ke(N)\cap F^{i-1}$. By assumption this means $N^{i-1}z=x-N^{i-1}y$ for some $z$, i.e. $x=N^{i-1}(z-y)$.
  
\begin{pro}\label{dannja} Suppose $k$ is finite, $Y$ is proper and satisfies (MW). Suppose in addition that the natural map $\iota:H_{rig}^*(Y)\to H^*$ is strict with respect to the canonical \v{C}ech filtrations, i.e. $\iota(F_C^rH_{rig}^*(Y))=\iota(H_{rig}^*(Y))\cap F_C^rH^*$ for all $r$. Then we have $\bi(N^r)=F_C^rH^*$ for all $r\ge 0$.   
\end{pro}

{\sc Proof:} By \ref{nillem} we only need to show $\ke(N)\cap F_C^rH^*= \ke(N)\cap\bi(N^r)$ for all $r$, but this is obviously implied by strictness of $\iota$ together with \ref{hrigana}.\\ 

\addtocounter{satz}{1}{{{}} \arabic{section}.\arabic{satz}} (1) The condition on $\iota$ in \ref{dannja} is formulated in terms of cohomology over the log base $T$ (using the double complex $B^{\bullet}$). But note that, given a lift of $Y$ to a semistable scheme $X$ over a finite totally ramified extension ${\mathcal O}_K$ of $W(k)$, the filtration $F_C^rH^*$ becomes the canonical \v{C}ech filtration on $H_{dR}^{*}(X_K)$ induced from $(C)_S$ (see section \ref{lifsec}), and the condition on $\iota$ becomes an entirely analytic condition not involving log structures.\\(2) Assuming (MW), the condition on $\iota$ holds if convoluting the filtrations $P^{can}_{\bullet}$ and $F_C^{\bullet}$ on $B^{\bullet}$ commutes with passing to cohomology; by \ref{fiwbmo} this is equivalent to: the convolution of the filtrations $\ke(N^{\bullet})$ and $F_C^{\bullet}H^*$ on $H^*$ is a weight filtration on $H^*$.\\
 
\addtocounter{satz}{1}{{{}} \arabic{section}.\arabic{satz}} Examples. (1) For curves we always have $\bi(N)=F_C^1H^1$, see \cite{coliov}.\\(2) Let $\Omega_K^{(d+1)}$ be Drinfel'd's $p$-adic symmetric space of dimension $d$ over a finite totally ramified extension $K$ of $K_0$, let $X_{\Gamma}$ be the quotient of $\Omega_K^{(d+1)}$ by a cocompact discrete torsionfree subgroup $\Gamma<PGL_{d+1}(K)$ and let $Y$ be its strictly semistable reduction. It has interesting cohomology only in degree $*=d$ where it satifies (MW) as was recently proven by de Shalit \cite{desh} and (independently) by Ito \cite{ito}. By \cite{hkstrat}, if $d$ is odd, we have $\bi(N^r)=F_C^rH^*$ for all $r$. But if $d$ is even we do {\it not always have $\bi(N^r)=F_C^rH^*$}. For example, let $d=2$. On $H^2_{dR}(X_{\Gamma})$ we have a covering filtration $(F^r_{\Gamma})_{r\ge 0}$ and in \cite{hkstrat} it is shown that it coincides with the \v{C}ech filtration $(F^r_{C})_{r\ge 0}$ on $H^2_{dR}(X_{\Gamma})$ (see section \ref{lifsec}), and moreover that $$F^j_{\Gamma}=\sum_i\ke(N^{i+1})\cap \bi(N^{i-2+2j})=\sum_i\ke(N^{i+1})\cap \bi(N^{i-3+2j})$$for the monodromy operator $N$ on $H^2_{dR}(X_{\Gamma})$ induced from $H^2_{crys}(Y/T)$ by means of $H^2_{dR}(X_{\Gamma})\cong H^2_{crys}(Y/T)\otimes_{W(k)}K=H^2\otimes_{K_0}K$. Observing $N^3=0$ we get $$F^1_{\Gamma}=\ke(N)+\bi(N)$$$$F^2_{\Gamma}=\ke(N)\cap\bi(N)+\bi(N^2)=\bi(N^2).$$ Now {\it assume} $\bi(N)=F_{C}^1H^2$ in $H^2=H^2_{crys}(Y/T)_{\mathbb{Q}}$. Then $\bi(N)=F^1_{\Gamma}$ in $H^2_{dR}(X_{\Gamma})$. Combining this with the above identities easily leads to $\ke(N)=\bi(N^2)$. Therefore $N$ induces an isomorphism $\bi(N)/\bi(N^2)\cong\bi(N^2)$. On the other hand, by the computations in \cite{ss} p.93 we have$$\dim_K(\bi(N^2))=\dim_K(F_{\Gamma}^2)=\dim_K(F^1_{\Gamma}/F_{\Gamma}^2)-1=\dim_K(\bi(N)/\bi(N^2))-1.$$Together this is a contradiction, disproving our assumption. 

\section{$N$ and singular cohomology}
\label{nundsing}

\addtocounter{satz}{1}{{{}} \arabic{section}.\arabic{satz}} In this section, $k$ is finite and $Y$ is proper. For $j\ge1$, irreducible components $M$ of $Y^j$ and irreducible components $N$ of $Y^{j+1}$ with $N\subset M$, denote by $$c^s_{M,N}:H_{crys}^s(M)_{\mathbb{Q}}\longrightarrow H_{crys}^s(N)_{\mathbb{Q}}$$the natural restriction maps. We say $Y$ is of weak Lefschetz type if for every $j\ge1$, every pair $(N,M)$ as above, the maps $c^s_{M,N}$ are isomorphisms if $s<(\dim N)=d-j$, and if the map $c^{d-j}_{M,N}$ is injective, with $\bi(c^{d-j}_{M,N})=h^{d-j}_N\subset H_{crys}^{d-j}(N)_{\mathbb{Q}}$ independent of $M$. For example, by \cite{km}, $Y$ is of weak Lefschetz type if it is projective and all embeddings $N\subset M$ as above are ample divisors on $M$.

\begin{lem}\label{weakleflem} Suppose $Y$ is of weak Lefschetz type, and that for each $i\ge1$, each component of $Y^i$ is geometrically connected. Let $s\le d-1$ and suppose $H^i(Y_{Zar},\mathbb{Q})=0$ for all $d-s>i>0$. Then the following sequence is exact:\begin{gather} H_{crys}^s(Y^1)_{\mathbb{Q}}\longrightarrow H_{crys}^s(Y^2)_{\mathbb{Q}}\longrightarrow\ldots\longrightarrow H_{crys}^s(Y^{d-s+1})_{\mathbb{Q}}\tag{$*$}\end{gather}
\end{lem}

{\sc Proof:} The geometrical connectedness of each component $L$ of $Y^i$ for each $i$ implies $H_{crys}^0(L)_{\mathbb{Q}}=K_0$ for each such $L$. Therefore the complex (with $H_{crys}^0(Y^1)_{\mathbb{Q}}$ in degree $0$)
\begin{gather}H_{crys}^0(Y^1)_{\mathbb{Q}}\longrightarrow H_{crys}^0(Y^2)_{\mathbb{Q}}\longrightarrow\ldots\longrightarrow H_{crys}^0(Y^{d+1})_{\mathbb{Q}}\longrightarrow0\tag{$**$}\end{gather} computes $H^*(Y_{Zar},K_0)=H^*(Y_{Zar},\mathbb{Q})\otimes K_0$, thus is exact by our hypothesis when truncated after the degree $(d-s)$ term. Now note that the weak Lefschetz assumption allows us to naturally identify into one single object $H^s$ all the following $K_0$-vector spaces: the cohomology groups $H_{crys}^s(L)_{\mathbb{Q}}$ for all components $L$ of $Y^j$ for all $1\le j\le d-s$, and the subspaces $h^{s}_L\subset H_{crys}^s(L)_{\mathbb{Q}}$ for all components $L$ of $Y^{d-s+1}$. We therefore obtain the sequence $(*)$ by tensoring $(**)$ over $K_0$ with $H^s$, truncating after the degree $(d-s)$ term and embedding the degree $(d-s)$ term into $H_{crys}^s(Y^{d-s+1})_{\mathbb{Q}}$. 

\begin{satz}\label{topvan} Suppose $Y$ satisfies (MW), is of weak Lefschetz type, and that for each $i\ge1$, each component of $Y^i$ is geometrically connected. Suppose $H^i(Y_{Zar},\mathbb{Q})=0$ for all $d>i>0$. Then $N=0$ on $H_{crys}^s(Y/T)_{\mathbb{Q}}$ for all $s\ne d=\dim(Y)$.
\end{satz}

{\sc Proof:} First note that we may assume $s<d$ since the assertion for $s>d$ is reduced to that for $2d-s$ using Poincar\'{e}-duality (which commutes with Frobenius, hence with weight filtrations, hence --- assuming (MW) --- with monodromy filtrations). By \ref{nillem} we need to show $\ke(N)\cap \bi(N)=(0)$ in $H_{crys}^s(Y/T)_{\mathbb{Q}}$. By \ref{hrigana} we can do this by proving $F^1_CH_{rig}^s(Y)=0$ (since we assume (MW)). We prove $F^t_CH_{rig}^s(Y)=0$ for $t\ge 1$ by descending induction on $t$. Proving $F^t_CH_{rig}^s(Y)=0$ means proving that $$\lambda_t:H^s(Y,P^{can}_1F_C^tA^{\bullet})\longrightarrow H^s(Y,P^{can}_1A^{\bullet})=H_{rig}^s(Y)$$is the zero map. In view of the exact sequence$$H^s(Y,P^{can}_1F_C^{t+1}A^{\bullet})\longrightarrow H^s(Y,P^{can}_1F_C^tA^{\bullet})\stackrel{\kappa}{\longrightarrow}H^s(Y,P^{can}_1\Gr_C^tA^{\bullet})$$and the vanishing of $\lambda_{t+1}$ by induction hypothesis, it is enough to show that for all $x\in H^s(Y,P^{can}_1F_C^tA^{\bullet})$ there exists a $x'\in H^s(Y,P^{can}_1F_C^tA^{\bullet})$ with $\kappa(x)=\kappa(x')$ and $\lambda_t(x')=0$. Consider the sequence \begin{gather}H^{s-1}(Y,P^{can}_1\Gr_C^{t-1}A^{\bullet})\stackrel{\alpha}{\longrightarrow}H^{s}(Y,P^{can}_1\Gr_C^tA^{\bullet})\stackrel{\beta}{\longrightarrow}H^{s+1}(Y,P^{can}_1\Gr_C^{t+1}A^{\bullet})\tag{$S$}\end{gather}where $\alpha$ and $\beta$ are the connecting maps in the obvious long exact cohomology sequences. Note $\Gr^j_CA^{\bullet}=A^{\bullet j}[-j]$ and $P^{can}_1A^{\bullet j}=\Gr_{j+1}{\bf C}\widetilde{\omega}^{\bullet}_Y[j+1]\cong {\bf C}\Omega^{\bullet}_{Y^{j+1}}$, together $P^{can}_1\Gr^j_CA^{\bullet}\cong {\bf C}\Omega^{\bullet}_{Y^{j+1}}[-j]$. Therefore $(S)$ becomes$$H_{crys}^{s-t}(Y^t)_{\mathbb{Q}}\longrightarrow H_{crys}^{s-t}(Y^{t+1})_{\mathbb{Q}}\longrightarrow H_{crys}^{s-t}(Y^{t+2})_{\mathbb{Q}},$$hence is exact by \ref{weakleflem}. Since $\kappa(x)\in\ke(\beta)$, we therefore find $y\in H^{s-1}(Y,P^{can}_1\Gr_C^{t-1}A^{\bullet})$ with $\alpha(y)=\kappa(x)$. The image $x'$ of $y$ under the composite $$H^{s-1}(Y,P^{can}_1\Gr_C^{t-1}A^{\bullet})\longrightarrow H^{s-1}(Y,P^{can}_1(A^{\bullet}/F_C^tA^{\bullet}))\longrightarrow H^{s}(Y,P^{can}_1F_C^tA^{\bullet})$$has the desired properties.\\

\addtocounter{satz}{1}{{{}} \arabic{section}.\arabic{satz}} Even if $H^i(Y_{Zar},\mathbb{Q})=0$ for all $i>0$, or even if $Y$ is the reduction of a semistable scheme $X$ over a finite totally ramified extension ${\mathcal O}_K$ of $W(k)$ whose generic fibre is contractible in the sense of Berkovich spaces (this hypothesis guarantees $H^i(Y_{Zar},\mathbb{Q})=0$ for all $d>i>0$ as required in \ref{topvan}), we can not expect the vanishing of $N$ on the middle degree cohomology $H_{crys}^d(Y/T)_{\mathbb{Q}}$. However, if we assume in addition that for every $j\ge1$, every irreducible component $M$ of $Y^j$ and every irreducible component $N$ of $Y^{j+1}$ with $N\subset M$, the restriction map$$H_{crys}^s(M)_{\mathbb{Q}}\longrightarrow H_{crys}^s(N)_{\mathbb{Q}}$$is an isomorphism if $s=d-j$, and is injective if $s=d-j+1$, with image independent on $M$, then we can argue as above to prove $N=0$ even on $H_{crys}^d(Y/T)_{\mathbb{Q}}$. For example we get $N=0$ on $H_{crys}^1(Y/T)_{\mathbb{Q}}$ in case $d=1$ (for this and a converse of it see also \cite{mokr} 5.6).\\

\addtocounter{satz}{1}{{{}} \arabic{section}.\arabic{satz}} If $Y_{\sigma}^i$ is geometrically connected for all $i\ge 1$, all $\sigma\in S_i$, then $H^s(Y_{Zar},K_0)=F_C^sH_{crys}^s(Y/T)_{\mathbb{Q}}$. Indeed, $F_C^sH_{crys}^s(Y/T)_{\mathbb{Q}}$ is the $E_{\infty}^{s,0}$-term of the spectral sequence \begin{gather}  E_1^{pq}=H^q(Y,B^{\bullet p})\Longrightarrow H^{p+q}(Y,B^{\bullet})=H^{p+q}_{crys}(Y/T)_{\mathbb{Q}}\tag*{$(C)_T$}\end{gather} from \arabic{cecmodn1}.\arabic{cecmodn2}. But $E_{\infty}^{s,0}=E_{2}^{s,0}$ is the $s$-th cohomology group of the complex$$H^0_{conv}(Y^{1}/T)\longrightarrow H^0_{conv}(Y^{2}/T)\longrightarrow H^0_{conv}(Y^{3}/T)\longrightarrow\ldots$$(with $H^0_{conv}(Y^{1}/T)$ in degree $0$). Our assumption implies $H^0_{conv}(Y^{i}_{\sigma}/T)=K_0$ for all $i\ge1$, all $\sigma\in S_i$, and therefore $E_{\infty}^{s,0}=H^s(Y_{Zar},K_0)$.\\Now suppose in addition $Y$ is of weak Lefschetz type and satisfies (MW). Suppose that we have $\bi(N^r)=F_C^rH^*$ (cf. section \ref{gleich}) and $N^r=0$ on $H_{crys}^r(Y/T)_{\mathbb{Q}}$ for all $0<r<d$. Then necessarily even $N=0$ on $H_{crys}^r(Y/T)_{\mathbb{Q}}$ for all $0<r<d$, as follows from \ref{topvan}.

\section{Liftings to mixed characteristic}
\label{lifsec}

\addtocounter{satz}{1}{{{}} \arabic{section}.\arabic{satz}}\newcounter{ssvgl1}\newcounter{ssvgl2}\setcounter{ssvgl1}{\value{section}}\setcounter{ssvgl2}{\value{satz}} Let $A$ be a complete discrete valuation ring which is a totally ramified finite extension of $W(k)$. Let $K=\quot(A)$ and fix a uniformizer $\pi$ in $A$. Let ${X}$ be a proper $\pi$-adic formal $\spf(A)$-scheme with strictly semistable reduction, i.e. Zariski locally it admits \'{e}tale maps to $\spf A<X_1,\ldots,X_{d+1}>/(X_1\ldots X_a-\pi)$ for some $1\le a\le d+1$. It is naturally a log smooth formal log-scheme over $$S=(\spf(A),(\mathbb{N}\longrightarrow A, 1\mapsto\pi)).$$ Taking reduction modulo $\pi$ we get a semistable $k$-log scheme $Y$ as in \arabic{dfsst1}.\arabic{dfsst2}. Let $X_K$ be the generic fibre of $X$ as a rigid analytic space. For $i\ge 0$ let $]Y^i[_{X}=\coprod_{r\in S_i}]Y_r^i[_{X}$, the direct sum of the tubes of the $Y_r^i$ in $X$. The covering $X_K=\cup_{j\in S_1}]Y^1_j[_X$ is an admissible open covering, it therefore gives rise to the \v{C}ech spectral sequence \begin{gather} E_1^{pq}=H_{dR}^q(]Y^{p+1}[_{X})\Longrightarrow H_{dR}^{p+q}(X_K).\tag*{$(C)_S$}\end{gather} On the other hand we have from \arabic{cecmodn1}.\arabic{cecmodn2} the spectral sequence \begin{gather}  E_1^{pq}=H^q(Y,B^{\bullet p})\Longrightarrow H^{p+q}(Y,B^{\bullet})=H^{p+q}_{crys}(Y/T)_{\mathbb{Q}}\tag*{$(C)_T$}.\end{gather}

\begin{pro} Depending on the choice of $\pi$ there is an isomorphism of spectral sequences$$(C)_S\cong (C)_T\otimes_{K_0}K.$$ 
\end{pro}

{\sc Proof:} Note $H^q(Y,B^{\bullet p})\cong H^q_{conv}(Y^{p+1}/T)$ by \ref{cocri} and $H^q_{crys}(Y/T)_{\mathbb{Q}}\cong H^q_{conv}(Y/T)$. Thus we get the isomorphism from \cite{hkstrat}. \\ 

\addtocounter{satz}{1}{{{}} \arabic{section}.\arabic{satz}}\newcounter{foprob1}\newcounter{foprob2}\setcounter{foprob1}{\value{section}}\setcounter{foprob2}{\value{satz}} By transport of structure the monodromy operator $N$ on $H_{crys}^*(Y/T)_{\mathbb{Q}}$ induces a monodromy operator $N$ on $H^*_{dR}(X_K)$ which does not depend on our choice of $\pi$, see \cite{hyoka} sect.5. Denote by $F_{C}^rH^*_{dR}(X_K)$ the filtration on $H^*_{dR}(X_K)$ induced by the spectral sequence $(C)_S$. Via $(C)_S\cong (C)_T\otimes_{K_0}K$ this is the filtration obtained by scalar extension from the filtration $F_C^{\bullet}H_{crys}^*(Y/T)_{\mathbb{Q}}$ on $H_{crys}^*(Y/T)_{\mathbb{Q}}$. Therefore we get from \ref{mokomptri} and \ref{inkl} the following theorem, which in particular gives an upper bound for the vanishing order of $N$ in terms of the rigid space ${X}_K$:

\begin{satz}\label{kinkl} There is a natural operator $N$ acting on $(C)_S$ inducing the monodromy operator $N$ on $H^*_{dR}(X_K)$. However, we have $N=0$ on all $E_1$-terms. In particular, $\bi(N^r)\subset N(F_{C}^{r-1}H^*_{dR}(X_K))\subset F_{C}^rH^*_{dR}(X_K)$ in $H^*_{dR}(X_K).$
\end{satz}

\addtocounter{satz}{1}{{{}} \arabic{section}.\arabic{satz}} We do not know if in general the residue map $H^{*}_{dR}({X}_K)\stackrel{Res}{\longrightarrow}F_C^i H^{*}_{dR}({X}_K))$ can be made explicit {\it without} involving the log basis $T$. However, for $i=*=d$ this map should be the following (generalizing that of \cite{coliov}): restrict a class in $H^{d}_{dR}({X}_K)$ to $H^{d}_{dR}(]Y^{d+1}[_{X})$; there choose a representing $d$-form, take its residue and view it as an element in $H^{0}_{dR}(]Y^{d+1}[_{X})$.\\We mention that also the tentative definition of $N$ given in \cite{alsh} for varieties $X_K$ uniformized by Drinfel'd's symmetric spaces is based on residue maps, and we expect that our description of $N$ can be used for a comparison with the $N$ from \cite{alsh}. 



\begin{thebibliography}{abcdefgh}  
\bibitem{alsh}{\it G. Alon and E. de Shalit}, Cohomology of discrete groups in harmonic cochains on buildings, Israel J. of Mathematics, {\bf 135} (2003), 355--377
\bibitem{bertat}{\it V. Berkovich}, An analog of Tate's conjecture over local and finitely generated fields, Internat. Math. Res. Notices {\bf 13} (2000), 665--680  
\bibitem{chia}{\it B. Chiarellotto}, Rigid cohomology and invariant cycles for a semistable log scheme, Duke Math. J. {\bf 97} (1999), no.1, 155--169 \bibitem{chile} {\it B. Le Stum and B. Chiarellotto}, Sur la puret\'{e} de la cohomologie cristalline, C. R. Acad. Sci. Paris, t 326, S\'{e}rie I, 961--963, 1998
\bibitem{coliov}{\it R. F. Coleman and A. Iovita}, The Frobenius and monodromy operators for curves and abelian varieties, Duke Math. J. {\bf 97} (1999), no.1, 171--215 
\bibitem{desh}{\it E. de Shalit}, The $p$-adic monodromy-weight conjecture for $p$-adically uniformized varieties, Compos. Math.  {\bf 141} (2005),  no. 1, 101--120. 
\bibitem{findag}\selectlanguage{german}{\it E. Gro"se-Kl\"onne}\selectlanguage{english}, Finiteness of de Rham cohomology in rigid analysis, Duke Math. J. {\bf 113} (2002), no.1, 57--91
\bibitem{colo}\selectlanguage{german}{\it E. Gro"se-Kl\"onne}\selectlanguage{english}, Compactification of log morphisms, Tohoku Math. J. {\bf 56} (2004), 79--104
\bibitem{hkstrat}\selectlanguage{german}{\it E. Gro"se-Kl\"onne}\selectlanguage{english}, Frobenius and Monodromy operators in rigid analysis, and Drinfel'd's symmetric space, Journal of Algebraic Geometry 14 (2005), 391--437
\bibitem{hy}{\it O. Hyodo}, On the de Rham-Witt complex attached to a semistable family, Comp. Math. {\bf 78} (1991), 241--260
\bibitem{hyoka} {\it O. Hyodo and K. Kato}, Semi-stable Reduction and Crystalline Cohomology with Logarithmic Poles, Asterisque 223, SMF, Paris (1994), 221-261
\bibitem{ilau}{\it L. Illusie}, Autour du th\'{e}or\`{e}me de monodromie locale, Astérisque No. 223, SMF, Paris (1994), 9--57
\bibitem{ito}{\it T. Ito}, Weight-Monodromy conjecture for $p$-adically uniformized varieties, Invent. Math. {\bf 159}  (2005),  no. 3, 607--656.  
\bibitem{fkato} {\it F. Kato}, Log smooth Deformation Theory. Tohoku Math. J. {\bf 48} (1996), 317--354
\bibitem{kalo} {\it K. Kato}, Logarithmic structures of Fontaine-Illusie, Algebraic Analysis, Geometry and Number Theory, J. Hopkins Univ. Press (1989), 191--224
\bibitem{km} {\it N. Katz, W. Messing}, Some consequences of the Riemann hypothesis for varieties over finite fields, Invent. Math. {\bf 23} (1974), 73--77
\bibitem{kiaub}{\it R. Kiehl}, Theorem A und Theorem B in der nichtarchimedischen Funktionentheorie, Invent. Math. {\bf 2} (1967),
256--273 
\bibitem{lestum}{\it B. Le Stum}, La structure de Hyodo-Kato pour les courbes,  Rend. Sem. Mat. Univ. Padova {\bf 94}  (1995), 279--301. 
\bibitem{mokr} {\it A. Mokrane}, La suite spectrale des poids en cohomologie de Hyodo-Kato, Duke Math. J. {\bf 72} (1993), 301-337
\bibitem{nakadeg}{\it Y. Nakkajima}, $p$-adic weight spectral sequences of log varieties, preprint 
\bibitem{ogid} {\it A. Ogus}, Logarithmic De Rham cohomology, preprint
\bibitem{ogcon}{\it A. Ogus}, $F$-crystals on schemes with constant log structure. Special issue in honour of Frans Oort, Compositio Math. {\bf 97} (1995), 187--225.
\bibitem{ss}{\it P. Schneider, U. Stuhler}, The cohomology of $p$-adic symmetric spaces, Inv. Math. {\bf 105}, 47--122 (1991)
\bibitem{shiho}{\it A. Shiho}, Crystalline fundamental groups. II. Log convergent and rigid cohomology, J. Math. Sci. Univ. Tokyo {\bf 9} (2002), 1--163.
\end{thebibliography}
\end{document}